\documentclass[11pt]{amsart}
\usepackage{amsmath,amssymb,mathrsfs,amscd}
\usepackage{graphicx}

\title{Explicit construction of new Moishezon twistor spaces}

\author{Nobuhiro Honda}
\thanks
{$^{\dag}$This work was partially supported by
Research Fellowships of the 
Japan Society for the Promotion
of Science for Young Scientists.
}
\date{}

\newcommand{\ol}{\overline}

\newcommand{\lra}{\longrightarrow}

\newcommand{\set}{\,|\,}
\newcommand{\proofend}{\hfill$\square$}

\newtheorem{prop}{Proposition}[section]

\newtheorem{thm}[prop]{Theorem}
\newtheorem{rmk}[prop]{Remark}
\newtheorem{cor}[prop]{Corollary}

\begin{document}

\begin{abstract}
In this paper we explicitly construct Moishezon twistor spaces on $n\mathbf{CP}^2$ for arbitrary $n\ge 2$ which admit a holomorphic $\mathbf C^*$-action.
When $n=2$, they coincide with Y. Poon's twistor spaces.
When $n=3$, they coincide with the ones studied by the author in \cite{Hon07-2}.
When $n\ge 4$, they are new twistor spaces, to the best of the author's knowledge.
By investigating the anticanonical system, we show that 
our twistor spaces are bimeromorphic to  conic bundles over certain rational surfaces.
The latter surfaces can be regarded as orbit spaces of the $\mathbf C^*$-action on the twistor spaces. Namely they are minitwistor spaces.
We explicitly determine their defining equations  in $\mathbf{CP}^4$.
It turns out that the structure of the minitwistor space is independent of $n$.
Further, we explicitly construct a $\mathbf{CP}^2$-bundle over the resolution of this surface, and provide an explicit defining equation of the conic bundles.
It shows that the number of irreducible components of the discriminant locus for the conic bundles increases as $n$ does.
Thus our twistor spaces have a lot of similarities with the famous LeBrun twistor spaces, where the minitwistor space $\mathbf{CP}^1\times\mathbf{CP}^1$ in LeBrun's case is replaced by our minitwistor spaces found in \cite{Hon-p-1}.
\end{abstract}

\maketitle

\bigskip\noindent
\section{ Introduction}
More than 15 years have passed since  C.\,LeBrun    \cite{LB91} discovered a series of  self-dual metrics  and their twistor spaces, on the connected sum of complex projective planes.
Basically they are obtained as a 1-dimensional reduction of the self-duality equation  \cite{AHS} for conformal classes, and
can be regarded as a hyperbolic version of gravitational  multi-instantons by G.\,Gibbons-S.\,Hawking \cite{GH78} and their twistor spaces  by N.\,Hitchin \cite{H79}.
%, and have played a central role in the study of compact twistor spaces.
Characteristic property of LeBrun's result is that it is completely explicit: for the twistor spaces, a bimeromorphic projective model is explicitly given by a defining equation and then bimeromorphic transformations are concretely given which produce actual twistor spaces.
Here  bimeromorphic transformations are essential, because compact twistor spaces cannot  be K\"ahler except two well known examples by a theorem of Hitchin \cite{Hi81}, and hence the projective model itself cannot be biholomorphic to the twistor spaces.

This explicitness  made it possible to handle LeBrun twistor spaces concretely and brought many knowledge about LeBrun twistor spaces and their small deformations.
For example, it is shown that the obstruction cohomology group  for deformations always vanishes
\cite{LB93}, and that the structure of LeBrun twistor spaces is stable under $\mathbf C^*$-equivariant deformations, at least if they possess only $\mathbf C^*$-symmetries  \cite{LB93,PP95}.
Non-general case when they possess  larger symmetries is also treated in
 \cite{Hon07-1} based on the explicit construction.
%Further,  if the LeBrun twistor spaces `degenerate' to admit 2-dimensional symmetries, 
%it was determined which subgroups in the isometry group yield new self-dual metrics with 1-dimensional symmetries \cite{Hon07-1}.
Furthermore,  many examples of Moishezon and non-Moishezon twistor spaces  were obtained  as small deformations of LeBrun twistor spaces \cite{Kr97, CK98, Hon01}.
Thus LeBrun twistor spaces have been the most important resource in the study of compact twistor spaces.
But unfortunately, once we shift to the deformed twistor spaces it is usually difficult to obtain their explicit construction, even when they can be shown to be Moishezon.

In this paper we would like to provide another resource by  explicitly constructing Moishezon twistor spaces on $n\mathbf{CP}^2$ for arbitrary $n\ge 2$, which admit  $\mathbf C^*$-symmetries.
When $n=2$ they coincide with the twistor spaces constructed by Y.\,Poon \cite{P86}.
When $n=3$ they coincide with the ones studied by the author in \cite{Hon07-2}.
If $n\ge 4$ they are new twistor spaces,  to the best of the author's knowledge.
Although they cannot be obtained as a small deformation of LeBrun twistor spaces for $n\ge 4$,
these twistor spaces have a number of common properties with LeBrun twistor spaces.
For example,  bimeromorphic projective models of the twistor spaces can be naturally realized as   conic bundles over certain rational surfaces.
Further, most significantly, the latter surfaces can be regarded as  minitwistor spaces in the sense of Hitchin \cite{H82-1,H82-2},  
whose structure was recently studied by the author in \cite{Hon-p-1}.
The structure of our minitwistor spaces is independent of $n$
(although they have a non-trivial moduli).
Instead, the structure of the discriminant locus depends on $n$.
%Thus our series of  twistor spaces can be viewed as a natural generalization in the case $n=2$ \cite{P86} and $n=3$ \cite{Hon07-2}.
Thus it may  be possible to say that  our twistor spaces are (non-trivial) `variants' or `cousins' of LeBrun twistor spaces.

In Section 2 we explain what kind of twistor spaces we are concerned.
We characterize our twistor spaces by the property that they contain a certain smooth rational surface $S$ (explicitly constructed as a blown-up of $\mathbf{CP}^1\times\mathbf{CP}^1$) as a member of the system $|(-1/2)K_Z|$.
When $n\ge 4$ this condition immediately shows that $ |(-1/2)K_Z|$ is a pencil, and it means that they are different from LeBrun twistor spaces or twistor spaces investigated by Campana-Kreu\ss ler \cite{CK98} of a degenerate form \cite{Hon03}.
In order to obtain more detailed structure of our twistor spaces, since the system $|(-1/2)K_Z|$ is only a pencil, we need to consider the next line bundle, the anticanonical line bundle.
In general it is not easy to investigate the anticanonical system especially for the case $n>4$ because some cohomology group obstructs.
But in the present case we find another route to show that it is  4-dimensional linear system (Prop.\,\ref{prop-ac1}).
Further, we can show that the image of the anticanonical map is always a quartic surface in $\mathbf{CP}^4$ whose defining equations can be explicitly written
(Prop.\,\ref{prop-quotient}).

In Section 3, we provide a natural realization of  bimeromorphic projective models of our twistor spaces, as  conic bundles over the minimal resolution of the image quartic surface
obtained in Section 2. 
This realization is an analogue of that of LeBrun twistor spaces \cite[Section 7]{LB91}
(see also H.\,Kurke's paper \cite{Ku92}).
%We note that the intrinsic meaning of this realization (as conic bundles) is due to Kreussler-Kurke \cite{} (as far as the author can see by checking bibliographies), and we rely on their interpretation.
We give an explicit defining equation of the  conic bundles (Theorem \ref{thm-1}).
Roughly speaking, the conic bundles are uniquely determined by a set of  $(n-2)$ anticanonical curves in the surface which have a unique node respectively.
For the explicit realization of the conic bundles we need to give an elimination of the indeterminacy locus for the anticanonical map. % (which will be summarized in the diagram \eqref{string}).
This step is again an analogue of the elimination of the base locus of the system $|(-1/2)K_Z|$ for the LeBrun twistor spaces.
However, since the base locus of our anticanonical system is more complicated, the present elimination  requires several steps.

In Section 4 we investigate the bimeromorphic map from our twistor spaces to the conic bundles given in Section 3 further, and decompose it into a succession of blowing-ups and blowing-downs.
In particular, we obtain an explicit operations for obtaining our twistor spaces, starting from the projective models (=\,the conic bundles) of Section 3.
This is still an analogue of the case for LeBrun, but our construction is again more complicated, partly because compared to the case of LeBrun there are more divisors in the projective models which do not exist in the actual twistor spaces, and, at the same time, some divisors are lacking in the projective models. 
These mean that we need more blowing-ups and blowing-downs.

Section 5 consists of 3 subsections.
In \S 5.1 we first show that our twistor spaces (studied in Sections 2--4) can be obtained as an equivariant small deformation of the twistor space of a Joyce metric on $n\mathbf{CP}^2$ of a particular kind, where the equivariancy is with respect to some $\mathbf C^*$-subgroup of $\mathbf C^*\times\mathbf C^*$ acting on the twistor spaces of Joyce metrics.
This guarantees the existence of our twistor spaces.
Next we see that the structure of our twistor spaces is stable under $\mathbf C^*$-equivariant small deformations.
In  \S 5.2, we compute the dimension of the moduli space of our twistor spaces.
The conclusion is it is $(3n-6)$-dimensional, which is exactly the same as that for general LeBrun twistor spaces.
We also remark that when $n\ge 4$ our twistor spaces cannot be obtained as a small deformation of LeBrun twistor spaces of any kinds.
In \S 5.3 we discuss a lot of similarities and differences between our twistor spaces and LeBrun twistor spaces.

\bigskip\noindent
{\bf Notations and Conventions.}
In our construction of twistor spaces, we take a number of blowing-ups and blowing-downs for (algebraic) 3-folds. 
To save  notations we adapt the following convention.
If $\mu:X\to Y$ is a bimeromorphic morphism between  complex varieties and $W$ is a complex subspace in $X$, we write $W$ for the image $\mu(W)$  if the restriction $\mu|_W$ is still bimeromorphic.
Similarly, if $V$ is a complex subspace of $Y$, we write $Y$ for the strict transformation of $Y$ into $X$.
%So $W\subset X$ is not isomorphic to the original $W\subset Y$ in general, and the same for the strict transforms.
If $D$ is a divisor on a complex manifold $X$, $[D]$ means the associated line bundle.
The dimension of the complete linear system $|D|$ means $\dim H^0(X,[D])-1$.
The base locus is denoted by ${\rm Bs}\,|D|$.
If a Lie group $G$ acts on $X$ by means of biholomorphic transformations and the divisor $D$ is $G$-invariant,
$G$ naturally acts on the vector space $H^0(X,[D])$.
Then $H^0(X,[D])^{G}$ means the subspace of $G$-invariant sections.
Further $|D|^G$ means its associated linear system.
(So all members of $|D|^G$ are $G$-invariant.)

\section{Twistor spaces with $\mathbf C^*$-action which have a particular invariant divisor}

Let $Z$ be a twistor space on $n\mathbf{CP}^2$.
It is known that there exists no divisors on general $Z$ if $n\ge 5$ and hence no meromorphic function exists for general $Z$.
So far most study on $Z$ are done for which the half-anticanonical system $|(-1/2)K_Z|$
is non-empty.
In this respect there is a fundamental result of  Pedersen-Poon \cite{PP94} saying that a real irreducible member $S\in |(-1/2)K_Z|$ is always a smooth rational surface and it is
biholomorphic to $2n$ points blown-up of $\mathbf{CP}^1\times\mathbf{CP}^1$.
Further, the blowing-down can be chosen in such a way that it preserves the real structure and, in that case the resulting real structure is necessarily given by 
\begin{align}\label{rs5}
\text{(complex conjugation)}\times\text{(anti-podal)}.
\end{align}
We write $\mathscr O(1,0)=p_1^*\mathscr O(1)$ and $\mathscr O(0,1)=p_2^*\mathscr O(1)$ where $p_i$ is the projection $\mathbf{CP}^1\times\mathbf{CP}^1\to\mathbf{CP}^1$ to the $i$-th factor.
(So the pencil $|\mathscr O(1,0)|$ has a circle's worth of real members and $|\mathscr O(0,1)|$ does not have real members.)

In general, the complex structure of $S$, namely a configuration of $2n$  points (to be blown-up) on  $\mathbf{CP}^1\times\mathbf{CP}^1$, has a strong effect on  algebro-geometric  structures of the twistor space $Z$ in which $S$ is contained.
For example if $n$ points among $2n$ are located on one and the same non-real curve $C_1$ in $|\mathscr O(1,0)|$ (or $|\mathscr O(0,1)|$), then $Z$ is necessarily so called a LeBrun twistor space \cite{LB91}.
In this case, the base locus of the system $ |(-1/2)K_Z|$ is precisely $C_1\cup\ol{C}_1\,(\subset S\subset Z$), where we are keeping the convention that the strict transformations are denoted by the same notation,
%($\ol{C}_1$ denotes the image of $C_1$ under the real structure), 
and after blowing-up $Z$ along $C_1\cup \ol{C}_1$, 
the meromorphic map associated to the system becomes a morphism whose image is a non-degenerate quadratic surface in $\mathbf{CP}^3=\mathbf P H^0((-1/2)K_Z)^{\vee}$.
This meromorphic map from $Z$ to the quadratic surface can  also be regarded as a quotient map of the natural $\mathbf C^*$-action on the LeBrun twistor spaces.
Further, every members of $|(-1/2)K_Z|$ are $\mathbf C^*$-invariant and $C_1\cup\ol{C}_1$ is exactly the 1-dimensional components of the $\mathbf C^*$-fixed locus.

\begin{figure}
\includegraphics{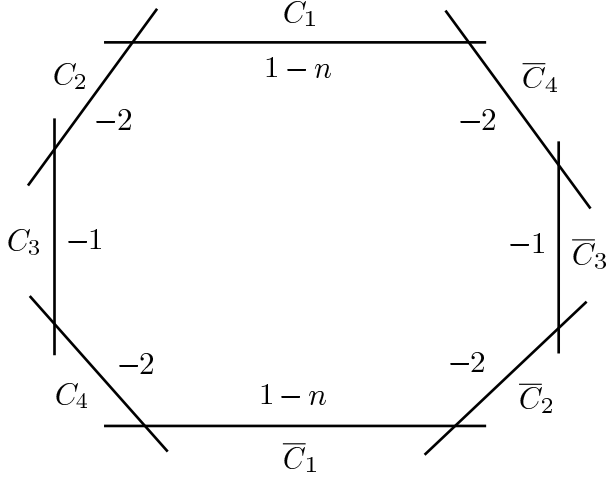}
\caption{The unique anticanonical curve $C$ on $S$}
\label{octagon}
\end{figure}

The starting point of the present investigation is to  consider a variant of  the above configuration of $2n$ points.
Namely we first choose a non-real curve $C_1\in|\mathscr O(1,0)|$ on $\mathbf{CP}^1\times\mathbf{CP}^1$ (as in the above LeBrun's case) and $(n-1)$ points $p_1,\cdots, p_{n-1}\in C_1$, and
let $S'$ be the blowing-up of $\mathbf{CP}^1\times\mathbf{CP}^1$ at $p_1,\cdots,p_{n-1},\ol{p}_1,\cdots,\ol{p}_{n-1}$, where $\ol{p}_i$ denotes the image of $p_i$ by the above real structure \eqref{rs5}.
Then $S'$ has an obvious non-trivial $\mathbf C^*$-action  which fixes every points of $C_1$ and $\ol{C}_1$.
Also there are $2(n-1)$ isolated fixed points on $S'$ which are on the exceptional curves of the blowing-up.
Among these $2(n-1)$ fixed points
we choose any two points which form a conjugate pair,  and let $S\to S'$ be the blowing-up at the points.
(We note that $n$ need to satisfy $n\ge 2$ for the construction works.)
$S$ has a lifted $\mathbf C^*$-action fixing $C_1$ and $\ol{C}_1$.
Further it is readily seen that if $n\ge 4$
the anticanonical system of $S$ consists of a unique member $C$ and it is a cycle of smooth rational curves consisting of 8 irreducible components, two of which are $C_1$ and $\ol{C}_1$.
We write 
\begin{equation}
C=\sum_{i=1}^4C_i+\sum_{i=1}^4\ol{C}_i,
\end{equation}
arranged as in Figure \ref{octagon}.
There, $C_3$ and $\ol{C}_3$ are the exceptional curves of the final blow-ups $S\to S'$.
It is also immediate to see from the above construction that their self-intersection numbers in $S$ satisfy
\begin{equation}
C_1^2=\ol{C}_1^2=1-n,\,\,C_2^2=\ol{C}_2^2=C_4^2=\ol{C}_4^2=-2, \,\,C_3^2=\ol{C}_3^2=-1.
\end{equation}

We are going to investigate  twistor spaces on $n\mathbf{CP}^2$ which have this rational surface $S$ as a divisor in $|(-1/2)K_Z|$. 
(The existence of such twistor spaces will be shown in the final section.)
As a preliminary we begin with the following

\begin{prop}\label{prop-S}
Let $S$ be the surface with $\mathbf C^*$-action (with $c_1^2=8-2n$) constructed above, and $C=\sum C_i+\sum \ol{C}_i$ the unique anticanonical curve of $S$.
Suppose $n\ge 4$.
Then (i) the fixed component of $|-2K_S|$ is $2C_1+C_2+C_4+2\ol{C}_1+\ol{C}_2+\ol{C}_4$,
(ii) the movable part of $|-2K_S|$ is free, 2-dimensional, and the image of the associated morphism $S\to\mathbf{CP}^2$ is a conic.
\end{prop}

\noindent Proof.
By computing intersection numbers, it is readily seen that the reducible  curve  in (i) is  contained in the fixed component of $|-2K_S|$.
Subtracting this from $-2K_S=2C$, we obtain the system $|C_2+2C_3+C_4+\ol{C}_2+2\ol{C}_3+\ol{C}_4|$.
Considering connected components, the latter system is generated by two systems  $|C_2+2C_3+C_4|$ and $|\ol{C}_2+2\ol{C}_3+\ol{C}_4|$.
Both of these are the pull-back of $|\mathscr O(0,1)|$ by the blowing-up $S\to\mathbf{CP}^1\times\mathbf{CP}^1$.
Hence their composite  is free, 2-dimensional and the image must be a conic.
\proofend

\bigskip
The structure of the half-anticanonical system $|(-1/2)K_Z|$ is as follows:

\begin{prop}\label{prop-fs}
Let   $Z$ be  a twistor space on $n\mathbf{CP}^2$, $n\ge 4$, equipped with a holomorphic $\mathbf C^*$-action compatible with the real structure, and suppose that there is a real  $\mathbf C^*$-invariant divisor\, $S\in |(-1/2)K_Z|$ which is equivariantly isomorphic to the complex surface in Prop.\,\ref{prop-S}.
Then the system $ |(-1/2)K_Z|$ satisfies the following.
(i) $\dim  |(-1/2)K_Z|=\dim  |(-1/2)K_Z|^{\mathbf C^*}=1$, and {\rm Bs}\,$|(-1/2)K_Z|=C$.
(ii) $ |(-1/2)K_Z|$ has precisely $4$ reducible members $S_i$ ($1\le i\le 4$), and they are   of the form $S_i=S_i^++S_i^-$, where $S_i^+$ and $S_i^-$ are mutually conjugate $\mathbf C^*$-invariant smooth divisors whose intersection numbers with twistor lines are one,
(iii) $L_i:=S_i^+\cap S_i^-$ $(1\le i\le 4)$ are\, $\mathbf C^*$-invariant twistor lines joining conjugate pairs of singular points of the cycle $C$.
\end{prop}

We omit a proof of this proposition since it is now standard and requires no new idea (cf.\,\cite[Prop.\,3.7]{Kr98} or \cite[Proof of Prop.\,1.2]{Hon99} for the proof of (ii)).
We distinguish irreducible components of  the 4 reducible members as follows: $S_1^+$ contains $C_1\cup C_2\cup C_3\cup C_4$, $S_2^+$ contains  $C_2\cup C_3\cup C_4\cup\ol{C}_1$, $S_3^+$ contains $C_3\cup C_4\cup\ol{C}_1\cup\ol{C}_2$ and $S_4^+$ contains $ C_4\cup\ol{C}_1\cup\ol{C}_2\cup\ol{C}_3$.
The $\mathbf C^*$-invariant cycle $C$ and these eight $\mathbf C^*$-invariant divisors $S_i^+$ and $S_i^-$ will repeatedly appear in our investigation of the twistor spaces.

For LeBrun twistor spaces, the system $|(-1/2)K_Z|$ is 3-dimensional and 
algebraic structures of the twistor spaces can be studied through its associated meromorphic map.
For our twistor spaces, since the system is only a pencil as in Prop.\,\ref{prop-fs}, we cannot go ahead if we consider the system only.
So we are going to study the next natural linear system, the anticanonical system.
The following proposition clarifies the structure of the anticanonical system of 
our twistor spaces, and plays a fundamental role throughout this paper.

\begin{prop}\label{prop-ac1}
Let $Z$ and  $S$ be as in Prop.\,\ref{prop-fs}.
Then  the following hold.
%(i) $\dim H^0((-1/2)K_Z)^{\mathbf C^*}=2$ for\, $n\ge 3$ and\, $H^0((-1/2)K_Z)^{\mathbf C^*}=H^0((-1/2)K_Z)$ for $n\ge 4$.
(i) $\dim |-K_Z|=\dim|-K_Z|^{\mathbf C^*}=4$.
(ii) The vector space $H^0(-K_Z)$ is generated by the image of a natural bilinear map
\begin{equation}\label{natmap}
H^0((-1/2)K_Z)\times H^0((-1/2)K_Z)\,\lra\, H^0(-K_Z)
\end{equation}
which generates a 3-dimensional linear subspace $V$ in $H^0(-K_Z)$, and 2 sections of \,$-K_Z$ defining the following 2 divisors 
\begin{equation}\label{ac1}
Y:=S_1^++S_2^++S_3^++S_4^-,\,\,\ol{Y}:=S_1^-+S_2^-+S_3^-+S_4^+.
\end{equation}
\end{prop}

\noindent Proof.
From the exact sequence 
\begin{equation}
0\,\lra (-1/2)K_Z\,\lra\,-K_Z\,\lra\,-2K_S\lra\,0,
\end{equation}
we obtain an exact sequence
\begin{equation}\label{es1}
0\,\lra H^0((-1/2)K_Z)\,\lra\,H^0(-K_Z)\,\lra\,H^0(-2K_S).
\end{equation}
To show surjectivity of the restriction map in \eqref{es1}, when $n=4$ we can use the Riemann-Roch formula and Hitchin's vanishing theorem \cite{Hi80} to deduce $H^1((-1/2)K_Z)=0$.
But the same calculation shows $H^1((-1/2)K_Z)\neq 0$ when $n\ge 5$.
This is the main difficulty when   investigating  algebraic structures of twistor spaces
on $n\mathbf{CP}^2$ in the case $n\ge 5$.
But in the present case we can proceed as follows.

We have $\dim H^0((-1/2)K_Z)=2$ by Prop.\,\ref{prop-fs} (i) and $\dim H^0(-2K_S)=3$ by Prop.\,\ref{prop-S} (ii).
Hence by \eqref{es1} we have $\dim H^0(-K_Z)\le 2+3=5$.
Since $\dim H^0((-1/2)K_Z)=2$, the subspace $V$ generated by the image of the bilinear map \eqref{natmap} is 3-dimensional.
To show $\dim H^0(-K_Z)=5$ it suffices to see that the 2 divisors $Y$ and $\ol{Y}$ in  \eqref{ac1} are actually anticanonical divisors on $Z$, and that 2 sections $y$ and $\ol{y}\in H^0(-K_Z)$ defining $Y$ and $\ol{Y}$  respectively satisfy $y\not\in V$, $\ol{y}\not\in V+\mathbf Cy$.
The former claim $Y,\ol{Y}\in|-K_Z|$ can be verified by computing the first Chern classes of the divisors $S_i^{\pm}$ (in $H^2(Z,\mathbf Z))$ explicitly as in \cite[Proof of Prop.\,1.2]{Hon99}.
So here we do not repeat the computations.
For the latter claims,  $y\not\in V$ is a direct consequence of  the fact that  $Y$ is not of the form $S+S'$ with $S,S'\in|(-1/2)K_Z|$,
which can again be verified by explicit forms of the first Chern classes of $S_i^{\pm}$.
%(Note that the divisors of the form $S_i^{\pm}+S_j^{\pm}$ $(1\le i,j\le 4$) is in $|(-1/2)K_Z|$ iff it is of the form $S_i^++S_i^-$ since otherwise 
In order to show $\ol{y}\not\in V+\mathbf Cy$, we first note that the base locus of the system $|V|$ is obviously the cycle $C$.
On the other hand from the explicit form \eqref{ac1} of $Y$ we obtain
  $Y|_S=\ol{C}_4+2C_1+3C_2+4C_3+3C_4+2\ol{C}_1+\ol{C}_2$.
But $\ol{Y}$ does not contain $C_3$.
These imply $\ol{y}\not\in V+\mathbf Cy$.
Thus we obtain  $\dim H^0(-K_Z)=5$.
We in particular obtain that the restriction map in \eqref{es1} is surjective for any $n\ge 4$.
Namely we have an exact sequence
\begin{equation}\label{es2}
0\,\lra H^0((-1/2)K_Z)\,\lra\,H^0(-K_Z)\,\lra\,H^0(-2K_S)\,\lra 0.
\end{equation}
To finish a proof of the proposition, it remains to show that 
$H^0(-K_Z)=H^0(-K_Z)^{\mathbf C^*}$.
By Prop.\,\ref{prop-fs} (i) and the equivariant exact sequence \eqref{es2}, it suffices to show that  the natural $\mathbf C^*$-action on $H^0(-2K_S)$ is trivial.
As seen in the proof of Prop.\,\ref{prop-S}, the movable part of $|-2K_S|$ is 
generated by two curves $C_2+2C_3+C_4$ and $\ol{C}_2+2\ol{C}_3+\ol{C}_4$ which are mutually linearly equivalent.
On the complex surface $S$, the line bundle $[C_2+2C_3+C_4]$ (and $[\ol{C}_2+2\ol{C}_3+\ol{C}_4]$) is isomorphic to the pullback of $|\mathscr O(0,1)|$ by the blowing-up $S\to\mathbf{CP}^1\times\mathbf{CP}^1$.
By our construction, $\mathbf C^*$ acts trivially on the second factor of $\mathbf{CP}^1\times\mathbf{CP}^1$.
Hence $\mathbf C^*$ acts trivially on $H^0(\mathscr O(0,1))$.
This implies that $\mathbf C^*$ acts trivially on $H^0(-2K_S)$, as required.
\proofend

\bigskip
As an easy consequence of Prop.\,\ref{prop-ac1} we obtain the following
\begin{cor}\label{cor-ac}
As generators of the system $|-K_Z|$ we can take  the following 5 divisors:
\begin{equation}
2S_1^++2S_1^-,\,2S_2^++2S_2^-,\,S_1^++S_1^-+S_2^++S_2^-\,\,\,\cdots { (\text{generator of  the system}}\,\,|V|),
\end{equation}
\begin{equation}
Y=S_1^++S_2^++S_3^++S_4^-,\,\,\ol{Y}=S_1^-+S_2^-+S_3^-+S_4^+.
\end{equation}
In particular, $|-K_Z|$ has no fixed component and its base locus is a curve
\begin{equation}
C-C_3-\ol{C}_3=(\ol{C}_4+C_1+C_2)+(C_4+\ol{C}_1+\ol{C}_2).
\end{equation}
\end{cor}

Next we study the meromorphic map associated to the anticanonical system.

\begin{prop}
\label{prop-quotient} 
Let $Z$ be as in Prop.\,\ref{prop-fs}, and $\Phi:Z\to\mathbf{CP}^4$ the meromorphic map associated to the system $|-K_Z|$.
Then we have the following.
(i)  The image \,$\mathscr T:=\Phi(Z)$ is a quartic surface defined by the two equations
\begin{equation}\label{T}
y_1y_2=y_0^2,\,\,\,y_3y_4=y_0\{y_1-\alpha y_2+(\alpha-1 )y_0\},
\end{equation}
where $(y_0,y_1,y_2,y_3,y_4)$  is a homogeneous coordinate on $\mathbf{CP}^4$, and $\alpha$ is  a real number satisfying $-1<\alpha<0$.
(ii) General fibers of \,$\Phi$ are the closures of  $\mathbf C^*$-orbits, and they are irreducible smooth  rational curves.
\end{prop}

\noindent Proof.
If $\sigma$ denotes the real structure on $Z$, the line bundle  $[S_i^-]$ ($1\le i\le 4)$ is isomorphic to the complex conjugation of $\sigma^*[S_i^+]$.
Let $e_i\in H^0([S_i^+])$ be a section which defines $S_i^+$, and $\ol{e}_i:=\ol{\sigma^* e_i}$ be a section of $[ S_i^-]$ defining  $S_i^-$.
Then $\{e_1\ol{e}_1,e_2\ol{e}_2\}$ is a basis of $H^0((-1/2)K_Z)\simeq\mathbf C^2$.
It is also a basis of the real part $H^0((-1/2)K_Z)^{\sigma}\simeq\mathbf R^2$.
Hence we can write
\begin{align}
e_3\ol{e}_3&=a\,e_1\ol{e}_1+b\,e_2\ol{e}_2,\label{fs1}\\
e_4\ol{e}_4&=c\,e_1\ol{e}_1+\alpha\,e_2\ol{e}_2\label{fs2},
\end{align}
for some $a,b,c,\alpha\in\mathbf R^{\times}$.
By multiplying  constants to $e_i$  we can suppose $a=c=1$ and $b=-1$.
Then 
noticing that $e_i\ol{e}_i\, (1\le i\le 4)$, considered as points of  $|(-1/2)K_Z|^{\sigma}\simeq\mathbf{RP}^1\simeq S^1$, are put in a linear order (i.e.\,clockwise or anti-clockwise order), 
we have $-1<\alpha<0$.
We set
\begin{equation}\label{012}
y_0=e_1e_2\ol{e}_1\ol{e}_2,\,\,y_1=(e_1\ol{e}_1)^2,\,\,y_2=(e_2\ol{e}_2)^2.
\end{equation}
These clearly form a basis of the 3-dimensional vector space $V\subset H^0(-K_Z)$.
Further we set
\begin{equation}
y_3=e_1e_2e_3\ol{e}_4,\,\,y_4=\ol{e}_1\ol{e}_2\ol{e}_3e_4,\,\,
\end{equation}
defining the 2 divisors $Y$ and $\ol{Y}$ respectively.
By Prop.\,\ref{prop-ac1}, $\{y_0,y_1,y_2,y_3,y_4\}$ forms a basis of $H^0(-K_Z)$.
In order to obtain a defining equation of the image $\mathscr T=\Phi(Z)$, we seek algebraic relations of this basis.
First we have an obvious relation
\begin{equation}\label{obv}
y_0^2=y_1y_2
\end{equation}
which immediately follows from \eqref{012}.
To obtain another relation, by \eqref{fs1} and \eqref{fs2}, we have
\begin{align}
e_3\ol{e}_3e_4\ol{e}_4&=(e_1\ol{e}_1-\,e_2\ol{e}_2)(e_1\ol{e}_1+\alpha\,e_2\ol{e}_2)\\
&=(e_1\ol{e}_1)^2-\alpha(e_2\ol{e}_2)^2+(\alpha-1)e_1\ol{e}_1e_2\ol{e}_2\\
&=y_1-\alpha y_2+(\alpha-1)y_0.
\end{align}
Hence we obtain
\begin{align}\label{nobv}
y_3y_4=e_1\ol{e}_1e_2\ol{e}_2e_3\ol{e}_3e_4\ol{e}_4
=y_0\{y_1-\alpha y_2+(\alpha-1)y_0\}.
\end{align}
Thus the image $\Phi(Z)\subset\mathbf{CP}^4$ is contained in the intersection of the two quadrics  \eqref{obv} and \eqref{nobv}.
To show that $\Phi$ maps surjectively to this quartic surface,
we note that there exists a diagram
\begin{equation}\label{cd1}
 \CD
Z@> \Phi>> \mathbf P H^0(-K_Z)^{\vee}\\
 @V\Psi VV @VV{\pi}V\\
\mathbf P H^0((-1/2)K_Z)^{\vee}@>{\iota}>>\mathbf P V^{\vee}\\
 \endCD
 \end{equation}
where $\Psi$ is the meromorphic map associated to the pencil $|(-1/2)K_Z|$,
$\pi$ is the projection induced by the inclusion $V\subset H^0(-K_Z)$,
and $\iota$ is the inclusion induced by the bilinear map \eqref{natmap} whose image is a conic defined by \eqref{obv}.
Further by the surjectivity of the restriction map in \eqref{es2}, 
the restriction of $\Phi$ on general fiber of $\Psi$ is precisely the map induced by the bi-anticanonical system on the fiber surface.
The image of the latter map is a conic (in the fiber plane of $\pi$) by Prop.\,\ref{prop-S} (ii).
Thus the image $\Phi(Z)$ is a (meromorphic) conic bundle over the conic (= the image of $\iota$).
On the other hand the intersection of \eqref{obv} and \eqref{nobv} also has an obvious structure of a (meromorphic) conic bundle structure over the same conic
(since \eqref{nobv} is quadratic).
This implies that $\Phi(Z)$ coincides with the intersection quartic surface.
Thus we obtain (i) of the proposition.

For (ii) we first note that the meromorphic map $\Phi$ is $\mathbf C^*$-equivariant (since it is associated to the anticanonical system) and the action on the target space is trivial by Prop.\,\ref{prop-ac1} (i).
Hence fibers of $\Phi$ are $\mathbf C^*$-invariant.
It remains to see the irreducibility and smoothness of general fibers.
Let $\Lambda$ be the image conic of $\iota$.
Then by the diagram \eqref{cd1} there is a natural rational map from $\mathscr T=\Phi(Z)$ to the conic $\Lambda$.
We still denote it by $\pi:\mathscr T\to\Lambda$.
Namely we have the following commutative diagram of meromorphic maps
\begin{equation}\label{cd2}
 \CD
Z@> \Phi>>\mathscr T\\
 @V\Psi VV @VV{\pi}V\\
\mathbf{CP}^1@>{\iota}>>\Lambda.\\
 \endCD
 \end{equation}
As above the restriction of $\Phi$ to a general fiber $\Psi^{-1}(\lambda)$ ($\lambda\in\mathbf{CP}^1$) is precisely the rational map associated to the bi-anticanonical system on the surface.
By Prop.\,\ref{prop-S}, after removing the fixed component, the bi-anticanonical system becomes free and a composite of two pencils whose general fibers are smooth rational curves.
Hence general fibers of $\Phi|_{\Psi^{-1}(\lambda)}:\Psi^{-1}(\lambda)\to\pi^{-1}(\lambda)\simeq\mathbf{CP}^1$ are irreducible and smooth.
This means that general fibers of $\Phi$ are smooth and irreducible.
Thus we have obtained all the claims of the proposition.
\proofend

\bigskip
Prop.\,\ref{prop-quotient} means that
the anticanonical map $\Phi:Z\to\mathscr T=\Phi(Z)$ is a (meromorphic) quotient map of the $\mathbf C^*$-action.
Thus our surface $\mathscr T$ is a parameter space of $\mathbf C^*$-orbits in $Z$.
Namely $\mathscr T$ is the {\em minitwistor space}\, in the sense of Hitchin \cite{H82-1,H82-2}.
Concerning the structure of $\mathscr T$ we have the following

\begin{prop}\label{prop-T} Let $\mathscr T$ be as in Prop.\,\ref{prop-quotient} and $\pi:\mathscr T\to \Lambda$ the projection to the conic as in the diagram \eqref{cd2}.
Then the following hold.
(i) The set of indeterminacy of\, $\pi$ consists of 2 points, and it coincides with the singular locus of $\mathscr T$.
Further both singularities are ordinary double points.
(ii) If $\tilde{\mathscr T}\to\mathscr T$ denotes the minimal resolution of the ODP's and $\tilde{\pi}:\tilde{\mathscr T}\to\Lambda$ denotes its composition with $\pi$, 
$\tilde{\pi}$ is a morphism whose general fibers are irreducible smooth rational curves.
Moreover, $\tilde{\pi}$ has precisely 4 singular fibers and all of them consist of two rational curves.
\end{prop}

Note that it follows from (ii) of the proposition that the resolved surface $\tilde{\mathscr T}$ is a rational surface satisfying $c_1^2=4$.
In the sequel we denote $\Gamma$ and $\ol{\Gamma}$ for the exceptional curves of the resolution $\tilde{\mathscr T}\to\mathscr T$.

\bigskip\noindent
Proof of Prop.\,\ref{prop-T}.
In the homogeneous coordinate $(y_0,y_1,y_2,y_3,y_4)$ in Prop.\,\ref{prop-quotient}, the projection $\pi$ is given by taking $(y_0,y_1,y_2)$.
Hence by the explicit equation \eqref{T}, $\pi$ is not defined only on the two points
$(0,0,0,1,0)\in\mathscr T$ and $(0,0,0,0,1)\in\mathscr T$.
Also it is elementary to see that these are ODP's of $\mathscr T$ and there are no other singularities of $\mathscr T$. Thus we obtain (i).

It is also elementary to see that $\pi$ becomes a morphism after taking the minimal resolution of the nodes.
The fibers of $\tilde{\pi}$ naturally correspond to those of $\pi$ and reducible ones of the latter are precisely over the intersection of the 2 conics $y_0^2=y_1y_2$ and $y_0\{y_1-\alpha y_2+(\alpha-1)y_0\}=0$.
This consists of 4 points $(0,1,0),(0,0,1),(1,1,1)$ and $(-\alpha,\alpha^2,1)$, and each fiber over there is two lines $y_3y_4=0$, as desired.
\proofend

\bigskip
We note that by the diagram \eqref{cd2} the conic $\Lambda=\{y_0^2=y_1y_2\}$ in $\mathbf{CP}^2=\mathbf PV^{\vee}$ (having $(y_0,y_1,y_2)$ as a homogeneous coordinate)  
is canonically identified with the parameter space of the pencil $|(-1/2)K_Z|$.
By the choice \eqref{012} of $y_0,y_1,y_2$, we have
\begin{align}
\Phi^{-1}(\{y_0=0\})=S_1^++S_1^-+S_2^++S_2^-,
\end{align}
\begin{align}
\Phi^{-1}(\{y_1=0\})=2S_1^++2S_1^-,\,\,
\Phi^{-1}(\{y_2=0\})=2S_2^++2S_2^-.
\end{align}
From these it follows that 
\begin{align}
\tilde{\pi}^{-1}(\{(0,0,1)\})=\Phi(S_1^+)\cup\Phi(S_1^-),\,\,
\tilde{\pi}^{-1}(\{(0,1,0)\})=\Phi(S_2^+)\cup\Phi(S_2^-).
\end{align}
Namely among the 4 critical values of  $\tilde{\pi}:\tilde{\mathscr T}\to\Lambda$,  the 2 points $(0,0,1)$ and $(0,1,0)$ correspond to the 2 reducible members $S_1^++S_1^-$ and $S_2^++S_2^-$  respectively.
On the other hand by \eqref{nobv} we have 
\begin{align}
\Phi^{-1}(\{y_1-\alpha y_2+(\alpha-1)y_0=0\})=S_3^++S_3^-+S_4^++S_4^-.
\end{align}
This means that the remaining 2 critical values $(1,1,1)$ and $(-\alpha,\alpha^2,1)$ correspond to the remaining 2 reducible members $S_3^++S_3^-$ and $S_4^++S_4^-$.

\begin{rmk}
{\em
The surface $\mathscr T$ was already investigated in \cite{Hon-p-1} as a minitwistor space of the twistor spaces with $\mathbf C^*$-action studied in \cite{Hon07-2}.
There, $\mathscr T$ is realized not as a quartic surface in $\mathbf{CP}^4$ but as a double covering of $\ol{\Sigma}_2$ (Hirzebruch surface of degree 2 whose $(-2)$ section is contracted), branched along a smooth elliptic (anticanonical) curve.
}
\end{rmk}

Since the $\mathbf C^*$-action on our twistor spaces is supposed to be compatible with the real structure, the corresponding self-dual structures on $n\mathbf{CP}^2$ carry $U(1)$-symmetry.
We finish this section by summarizing properties of this $U(1)$-action.

\begin{prop}\label{prop-U(1)}\label{prop-u1}
The $U(1)$-action on $n\mathbf{CP}^2$ induced from the $\mathbf C^*$-action on our twistor spaces satisfies the following.
(i) $U(1)$-fixed locus consists of one sphere and isolated $n$ points.
(ii) There exists a unique $U(1)$-invariant sphere along which the subgroup $\{\pm1\}\subset U(1)$ acts trivially.
Further, this sphere contains precisely 2 isolated $U(1)$ fixed points.
(iii) The 2 fixed points in (ii) can be characterized by the property that the twistor lines over the points are not fixed by the $\mathbf C^*$-action.
(Namely, the points on the twistor lines over other $(n-2)$ isolated fixed points are fixed.)
(iv) The $U(1)$-action is free outside the fixed locus and the $U(1)$-invariant sphere in (ii).
\end{prop}

This can be readily obtained by making use of our explicit $\mathbf C^*$-action on the real invariant divisor $S\in|(-1/2)K_Z|$ and the twistor fibration $Z\to n\mathbf{CP}^2$.
Note that our $U(1)$-action on $n\mathbf{CP}^2$ is unique since the $\mathbf C^*$-action on the invariant divisor $S$ is unique up to diffeomorphisms.
Of course, the  $U(1)$-fixed sphere in (i) is the image of the $\mathbf C^*$-fixed rational curve $C_1$ (and $\ol{C}_1$).
(ii) implies that our $U(1)$-action on $n\mathbf{CP}^2$ is not {\em semi-free},
but (iv) implies that it is almost semi-free.
The $U(1)$-invariant sphere having the isotropy subgroup $\{\pm 1\}$ is the image of the exceptional curves of the final blowing-up $S\to S'$ for obtaining $S$.
As for (iii) the two $\mathbf C^*$-invariant twistor lines which are not fixed are exactly $L_3=S_3^+\cap S_3^-$ and $L_4=S_4^+\cap S_4^-$.
The other invariant twistor lines $L_1=S_1^+\cap S_1^-$ and  $L_2=S_2^+\cap S_2^-$ are  over the $U(1)$-fixed sphere.
The remaining $(n-2)$ fixed   twistor lines go through isolated $\mathbf C^*$-fixed points on the divisor $S$.

We note that for the  $U(1)$-action of LeBrun metrics, (i) holds, but (ii)  does not hold.
Namely the action is semi-free, and the twistor lines over isolated fixed points are always $\mathbf C^*$-fixed.

\section{Defining equations of  projective models as conic bundles}

In the last section we showed that the anticanonical system of our twistor spaces gives a  meromorphic map $\Phi:Z\to\mathscr T$ and the image surface $\mathscr T$ can be viewed as a minitwistor space whose natural defining equations in $\mathbf{CP}^4$ can be explicitly written down.
In this section we investigate the map $\Phi$ more in detail and show that 
there is a natural bimeromorphic map from our twistor space to a certain conic bundle on the resolved minitwistor space $\tilde{\mathscr T}$ (cf.\,Prop.\,\ref{prop-T}).
Further, we explicitly construct a $\mathbf{CP}^2$-bundle over $\tilde{\mathscr T}$ in which the conic bundle is  embedded, and also give the defining equations of the conic bundle in the $\mathbf{CP}^2$-bundle.
Basically these are accomplished by eliminating the indeterminacy of the anticanonical map $\Phi:Z\to\mathscr T$.

\begin{figure}
\includegraphics{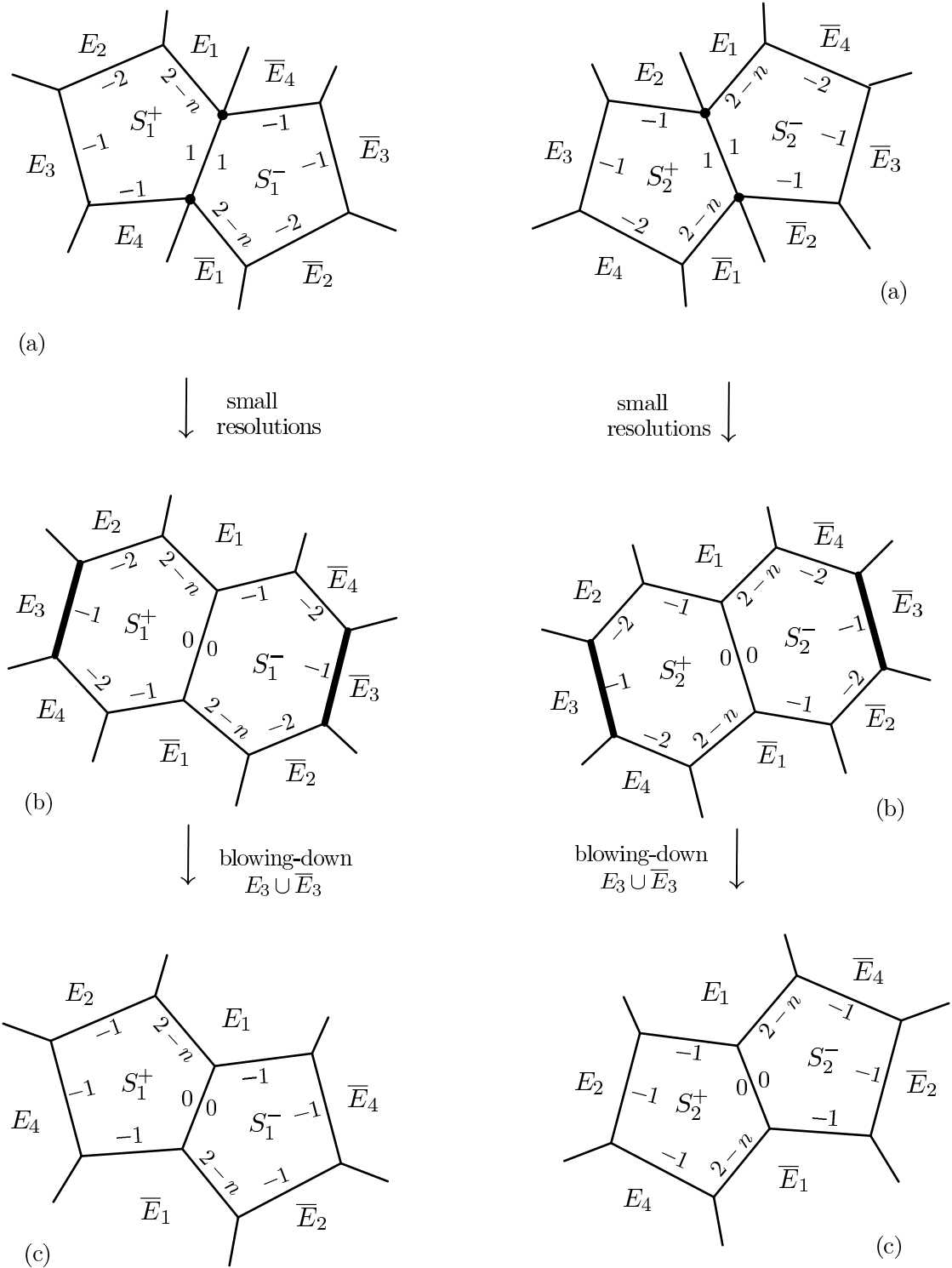}
\caption{operations for $S_1^{\pm}$ and $S_2^{\pm}$}
\label{fig1_and_2}
\end{figure}

\begin{figure}
\includegraphics{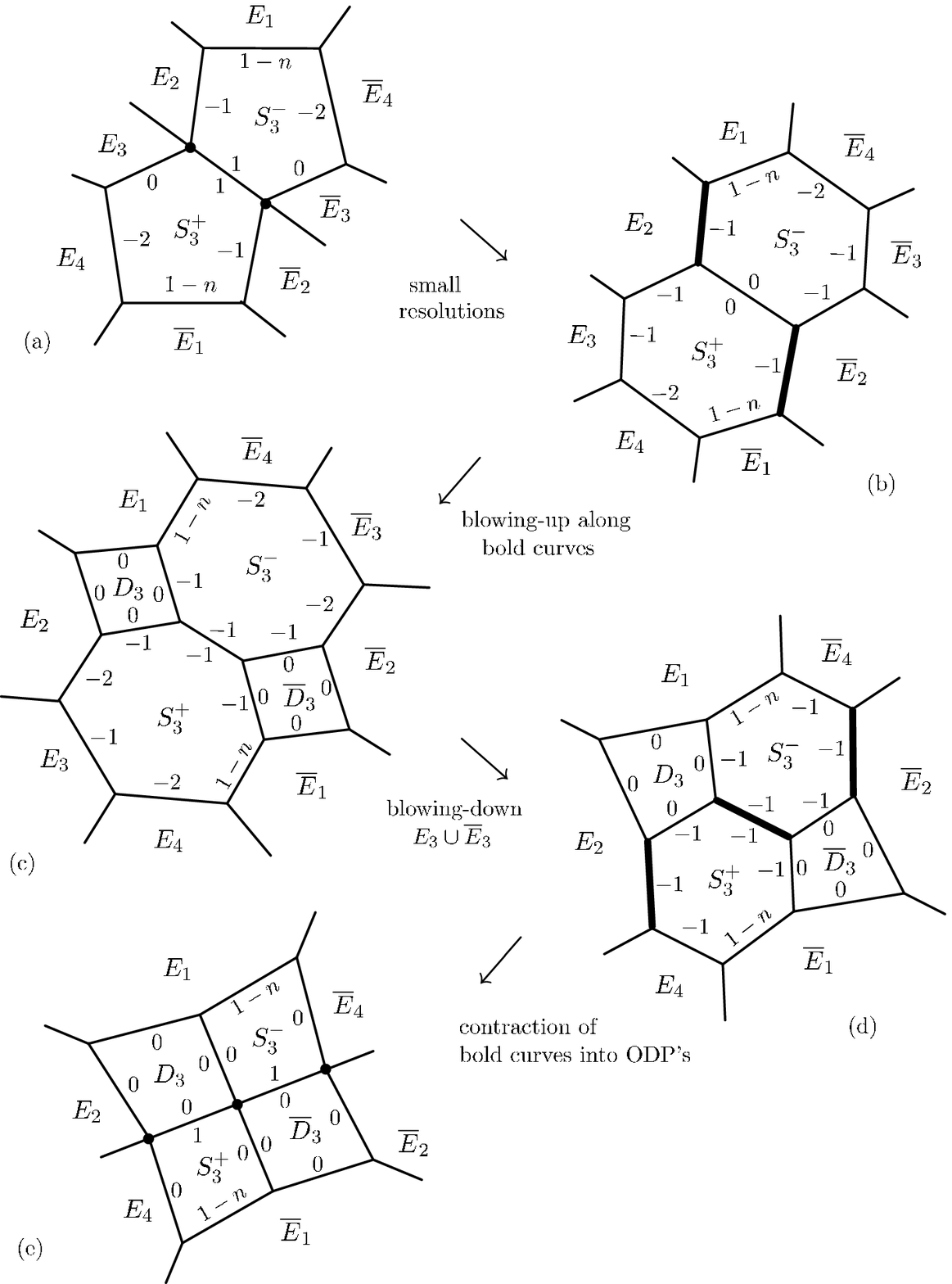}
\caption{operations for $S_3^{\pm}$}
\label{fig3}
\end{figure}

\begin{figure}
\includegraphics{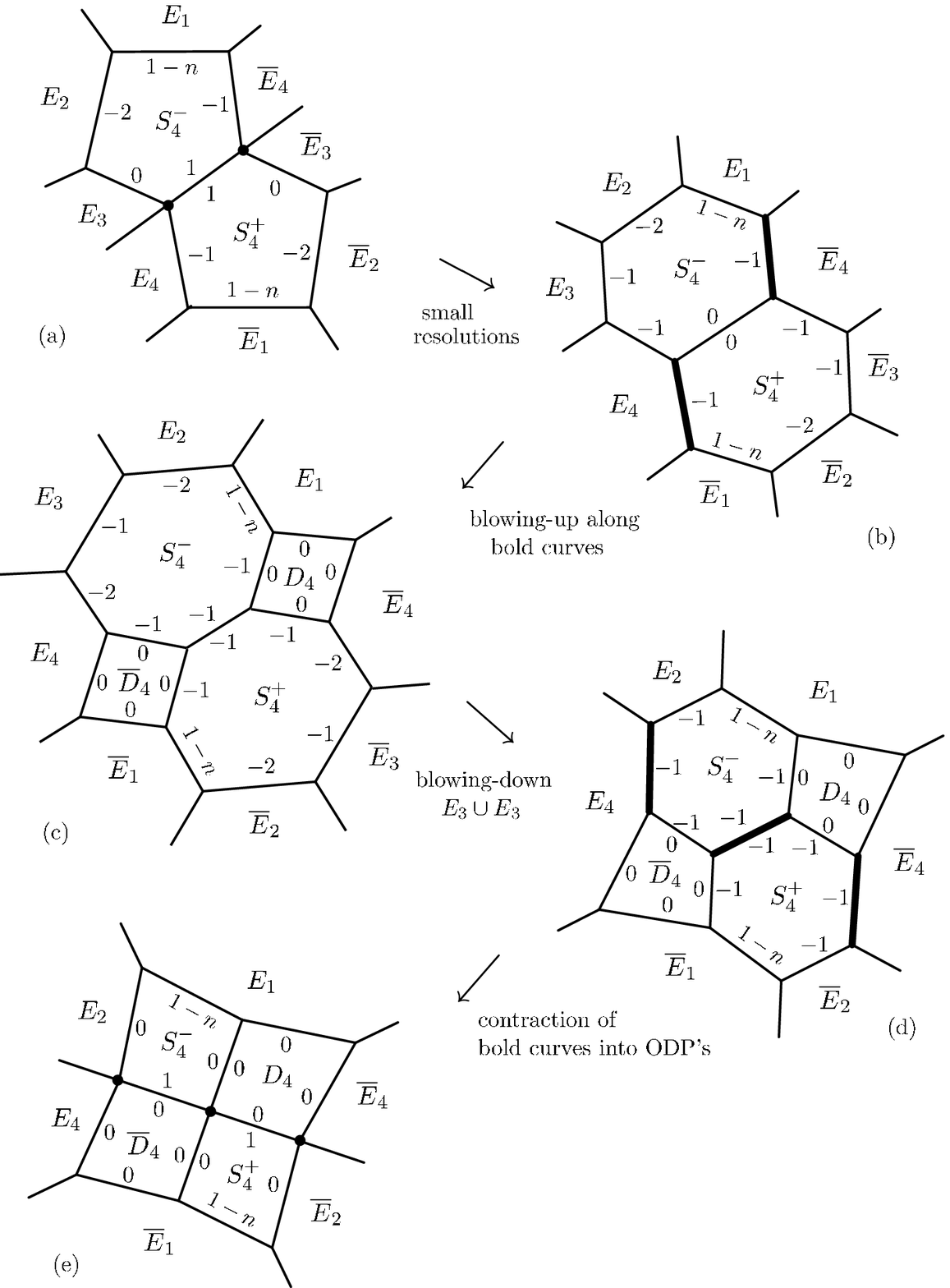}
\caption{operations for $S_4^{\pm}$}
\label{fig4}
\end{figure}

We are now going to eliminate the indeterminacy locus of $\Phi$ by a succession of blowing-ups, where we use Cor.\,\ref{cor-ac} to know where we have to blow-up.
As a first step let $Z_1\to Z$ be the blowing-up along the cycle $C$, and $E_i$ and $\ol{E}_i$ $(1\le i\le 4)$ the exceptional divisors over $C_i$ and $\ol{C}_i$ respectively.
Since $C$ has 8 nodes, $Z_1$ has 8 ordinary double points.
Each ODP is on the intersection of two exceptional divisors.
Noting that $\Phi$ maps the components $C_3$ and $\ol{C}_3$  to the 2 nodes
$\mathscr T$, $\Phi$ is naturally lifted to a (still non-holomorphic) map $\Phi_1:Z_1\to\tilde{\mathscr T}$ in such a way that the following diagram is commutative:
\begin{equation}\label{cd3}
 \CD
Z_1@>>>Z\\
 @V\Phi_1 VV @VV{\Phi}V\\
\tilde{\mathscr T}@>>>\mathscr T.\\
 \endCD
 \end{equation}
 Further, $\Phi_1$ satisfies
$\Phi_1(E_3)=\Gamma,\Phi_1(\ol{E}_3)=\ol{\Gamma}$,
where $\Gamma$ and $\ol{\Gamma}$ are the exceptional curves of the resolution $\tilde{\mathscr T}\to\mathscr T$ as before.
By construction, fibers of  the composition map $Z_1\to\tilde{\mathscr T}\to\Lambda\simeq\mathbf{CP}^1$ 
are the strict transforms of members of $|(-1/2)K_Z|$.
As in Prop.\,\ref{prop-fs} there are 4 real reducible members $S_i^++S_i^-$ ($1\le i\le 4$) of $|(-1/2)K_Z|$.
Correspondingly there are 4 real points of $\Lambda$ whose inverse images under the above map 
$Z_1\to\Lambda$ are the strict transforms of the 4 reducible members.
These are exactly the points where the morphism $\tilde{\mathscr T}\to\Lambda$ has reducible fibers.
The 4 reducible fibers of the map $Z_1\to\Lambda$ are illustrated in (a) of Figures \ref{fig1_and_2},\,\ref{fig3},\,\ref{fig4},  where the dotted points represent the ODP's of $Z_1$.
The points are precisely the points where 4 faces (representing $\mathbf C^*$-invariant divisors) meet.
%(These corresponds to the 4 reducible fibers of $\tilde{\pi}$.)

As a second step for elimination, we  take  small resolutions of all these ODP's of $Z_1$.
Of course there are 2 choices for each ODP's and we distinguish them by specifying which pair of divisors is blown-up at the shared ODP.
We choose the following small resolutions:
\begin{itemize}
\item
On the 2 ODP's on $L_1$, the pairs $\{S_1^+,\ol{E}_1\}$ and $\{S_1^-,E_1\}$ are blown-up.
\item
On the 2 ODP's on $L_2$, the pairs $\{S_2^+,E_1\}$ and $\{S_2^-,\ol{E}_1\}$ are blown-up.
\item
On the 2 ODP's on $L_3$, the pairs $\{S_3^+,E_2\}$ and $\{S_3^-,\ol{E}_2\}$ are blown-up.
\item
On  the 2 ODP's on $L_4$, the pairs $\{S_4^+,\ol{E}_4\}$ and $\{S_4^-,E_4\}$ are blown-up.
\end{itemize}
(See (a) $\to$ (b) of Figures \ref{fig1_and_2},\,\ref{fig3},\,\ref{fig4}.)
We denote $Z_2$ for the resulting   non-singular 3-fold.
Note that the small resolution $Z_2\to Z_1$ preserves the real structure.
We denote the composition map $Z_2\to Z_1\to\tilde{\mathscr T}$ by $\Phi_2$.
$\Phi_2$ is still non-holomorphic.
We denote the composition bimeromorphic morphism $Z_2\to Z_1\to Z$ by $\mu_2$.

Since we know explicit generating divisors of $|-K_Z|$ as in Cor.\,\ref{cor-ac} and the way how they intersect, we can 
chase the changes of the base locus of the system under  $\mu_2:Z_2\to Z$.
Namely, we have
\begin{align}\label{bs2}
{\rm Bs}\,
\left|\mu_2^*(-K_Z)-\sum_{i=1,2,4}(E_i+\ol{E}_i)\right|=
(S_3^-\cap E_2)\cup( S_3^+\cap\ol{E}_2)\cup(S_4^-\cap \ol{E}_4)\cup(S_4^+\cap E_4),
\end{align}
where each intersections on the right are smooth $\mathbf C^*$-invariant rational curves.
(In (b) of Figures \ref{fig3} and \ref{fig4}, these are illustrated by bold lines.)
So as a third step for elimination let $Z_3\to Z_2$ be the blowing-up along these 4 curves,
and $D_3,\ol{D}_3,D_4,\ol{D}_4$ the exceptional divisors respectively, named in the order of \eqref{bs2}.
Then it can also be seen that the system becomes free and hence the composite map $\Phi_3:Z_3\to( Z_2\to Z_1\to)\, \tilde{\mathscr T}$ is  a morphism.
Instead, the 2 fibers of the morphism $Z_3\to(\tilde{\mathscr T}\to)\,\Lambda$ containing $S_3^++S_3^-$ and $S_4^++S_4^-$ consist of 4 components
$S_3^++S_3^-+D_3+\ol{D}_3$ and $S_4^++S_4^-+D_4+\ol{D}_4$ respectively.
(See Figures \ref{fig3},\,\ref{fig4} (c).)

Thus we have obtained a sequence of blow-ups which eliminates the indeterminacy locus of the meromorphic map $\Phi$:
\begin{equation}\label{string}
 \CD
Z_3@>>> Z_2@>>>  Z_1@>>> Z\\
 @V{\Phi_3}VV      @V{\Phi_2}VV         @V{\Phi_1}VV @VV{\Phi}V\\
  \tilde{\mathscr T}        @=             \tilde{\mathscr T} @=\tilde{\mathscr T}@>>>\mathscr T,\\
 \endCD
 \end{equation}
 where among vertical arrows  only $\Phi_3$ is a morphism.
We note that all centers of blow-ups are $\mathbf C^*$-invariant and hence the whole of \eqref{string} preserves $\mathbf C^*$-actions.
Also all the blow-ups preserve the real structures.
Since $\mathbf C^*$ acts trivially on $\tilde{\mathscr T}$,
 all fibers of the morphism $\Phi_3$ are $\mathbf C^*$-invariant, and they  are generically smooth and irreducible since it is already true for the original map $\Phi$ as in Prop.\,\ref{prop-quotient} (ii).
Further,  $E_1$ and $\ol{E}_1$ in $Z_3$ are {\em sections}\, of $\Phi_3$.
In particular, they are biholomorphic to the surface $\tilde{\mathscr T}$.

\begin{figure}
\includegraphics{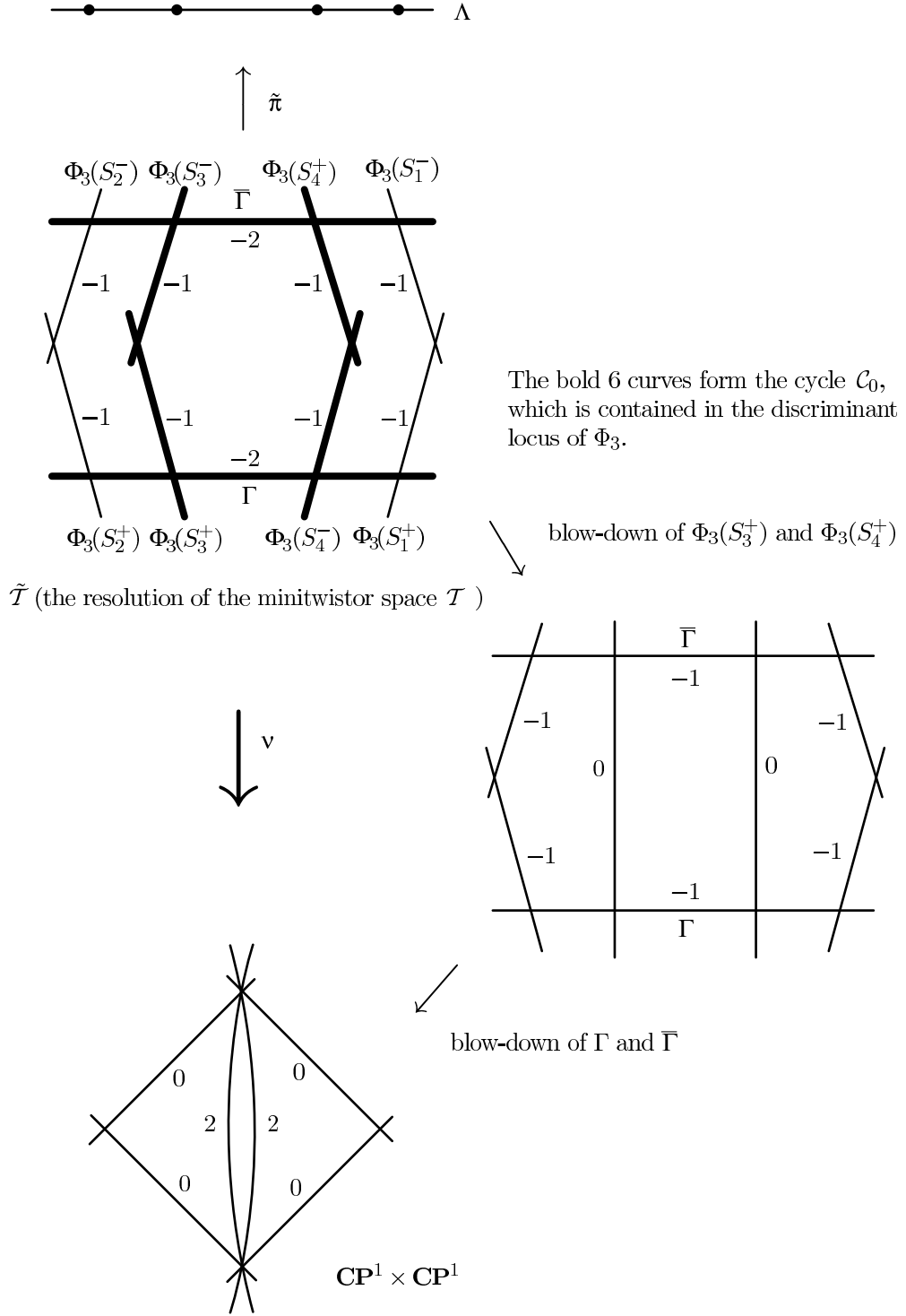}
\caption{structure of the minitwistor space $\mathscr T$ and illustration of $\nu$}
\label{figT}
\end{figure}

Next in order to express the normal bundles of $E_1$ and $\ol{E}_1$ in $Z_3$ in simple forms, we first realize
the surface $\tilde{\mathscr T}$ as  a blown-up of $\mathbf{CP}^1\times\mathbf{CP}^1$.
For this, we again recall that the conic bundle map $\tilde{\pi}:\tilde{\mathscr T}\to\Lambda$ has 4 reducible fibers and that they are precisely the images of the reducible members of the system $|(-1/2)K_Z|$ (cf.\,Prop.\,\ref{prop-T} and the explanation following its proof).
We first blow-down two of the irreducible components of reducible fibers of $\tilde{\pi}$:
explicitly, we blow-down the components $\Phi(S_3^+)$ and $\Phi(S_4^+)$.
(If we use  the morphism $\Phi_3$ instead, these are equal to $\Phi_3(D_3)$ and $\Phi_3(D_4)$ respectively.)
Consequently the curves $\Gamma$ and $\ol{\Gamma}$ become $(-1)$-curves.
So we blow-down these two $(-1)$-curves.
Then the resulting surface is  $\mathbf{CP}^1\times\mathbf{CP}^1$.
We denote the composition morphism by 
\begin{equation}\label{nu}
\nu:\tilde{\mathscr T}\to\mathbf{CP}^1\times\mathbf{CP}^1.
\end{equation}
(See Figure \ref{figT}.)
Note that $\nu$ never preserves the real structure since it blows down $(-1)$-curves contained in real fibers of $\tilde{\pi}$.
We distinguish two factors of $\mathbf{CP}^1\times\mathbf{CP}^1$ by declaring that the image curve $\nu(\Phi(S_1^+))$ is a fiber of the second projection.
(Then the curve $\nu(\Phi(S_1^-))$ is a fiber of the first projection.)
Then by our explicit way of blowing-ups  $\mu_3:Z_3\to Z$, we can verify that the normal bundle of $E_1$ in $Z_3$ satisfies
\begin{align}\label{nb-1}
N_{E_1/Z_3}\simeq\nu^*\mathscr O(-1,2-n),
\end{align}
where we are using the isomorphism $E_1\simeq \tilde{\mathscr T}$ induced by $\Phi_3$.
We denote the line bundle on $\tilde{\mathscr T}$ on the right side by $\mathscr N$.
Then by reality, we have
\begin{align}\label{nb-2}
N_{\ol{E}_1/Z_3}\simeq\ol{\sigma^*\mathscr N}=:\mathscr N',
\end{align}
where $\sigma$ denotes the natural real structure on $\tilde{\mathscr T}$ induced from that on the twistor space $Z$, and again we are using the isomorphism $\ol{E}_1\to\tilde{\mathscr T}$ induced by $\Phi_3$.

Next we embed our threefold $Z_3$ into a $\mathbf{CP}^2$-bundle over $\tilde{\mathscr T}$ as a conic bundle,  up to contractions of some divisors and rational curves.
For this, important thing is the discriminant locus of the morphism $\Phi_3:Z_3\to\tilde{\mathscr T}$.
By our explicit way of blowing-ups, we can find 6 smooth rational curves in $\tilde{\mathscr T}$ such that  their inverse images split into 2 or 3 irreducible components:
two of them are the $(-2)$-curves $\Gamma$ and $\ol{\Gamma}$.
In fact, for these  curves, we have
\begin{align}\label{dl3}
\Phi_3^{-1}(\Gamma)=E_2+E_3+E_4, \,\,
\Phi_3^{-1}(\ol{\Gamma})=\ol{E}_2+\ol{E}_3+\ol{E}_4
\end{align}
and each components are mapped surjectively to the curves.
Thus $\Gamma$ and $\ol{\Gamma}$ are contained in the discriminant locus of $\Phi_3$.
The other 4 curves are the (reducible) fibers of $\tilde{\pi}:\tilde{\mathscr T}\to\Lambda$ over the 2 points $(1,1,1)$ and $(-\alpha,\alpha^2,1)\in\Lambda$.
In fact, the inverse images of the irreducible components of these 2 fibers are 
$S_3^++D_3$, $S_3^-+\ol{D}_3$, $S_4^-+\ol{D}_4$ and $S_4^++D_4$, where $D_3,D_4,\ol{D}_3$ and $\ol{D}_4$ are, as before, concretely named exceptional divisors of the final blowing-up $Z_3\to Z_2$.
All of these components are mapped surjectively to (one of) the 4 curves.
This means that  the 4 curves (in $\tilde{\mathscr T}$) are also irreducible components of the discriminant locus of $\Phi_3$.
So we have obtained 6 irreducible components of the discriminant locus of $\Phi_3$.
Let $\mathscr C_0$ be the sum of these 6 curves.
It is clearly  a cycle of rational curves. (See Figure \ref{figT}.)
Then it is elementary to see that $\mathscr C_0$ is a real anticanonical curve of $\tilde{\mathscr T}$.

With these preliminary construction we have the following result which explicitly gives a projective models of our twistor spaces:

\begin{thm}\label{thm-1}
Let $Z$ be a twistor space on $n\mathbf{CP}^2$ with $\mathbf C^*$-action which has the complex surface $S$ (given in Section 2) as a real $\mathbf C^*$-invariant divisor in the system $|(-1/2)K_Z|$ as in Prop.\,\ref{prop-fs}.
Let $\mathscr T$ be the image surface of the anticanonical map $\Phi$ as in Prop.\,\ref{prop-quotient}, $\tilde{\mathscr T}\to\mathscr T$ the minimal resolution as in Prop.\,\ref{prop-T}, and $\nu:\tilde{\mathscr T}\to\mathbf{CP}^1\times\mathbf{CP}^1$ the blowing-down obtained in \eqref{nu}.
Then the following hold.
(i) $Z$ is bimeromorphic to a conic bundle 
\begin{align}\label{c-bdle}
xy=t^2P_0P_5P_6\cdots P_{n+2},
\end{align}
defined in the $\mathbf{CP}^2$-bundle $\mathbf P(\mathscr N^{\vee}\oplus\mathscr N'^{\vee}\oplus
 \mathscr O)\to\tilde{\mathscr T}$,
where 
\begin{equation}
x\in\mathscr N^{\vee}=\nu^*\mathscr O(1,n-2),\,\,y\in \mathscr N'^{\vee}=\ol{\sigma^*\mathscr N^{\vee}},\,\,t\in\mathscr O,
\end{equation}
and $P_0$ and $P_j$ ($5\le j\le n+2$) are non-zero sections of the anticanonical bundle of $\tilde{\mathscr T}$.
(ii) The curve $\{P_0=0\}$ coincides with the curve $\mathscr C_0$ introduced above.
(iii)  The anticanonical curves $\mathscr C_j:=\{P_j=0\}$ ($5\le j\le n+2$) are real, irreducible and have a unique node respectively.
\end{thm}

\begin{figure}
\includegraphics{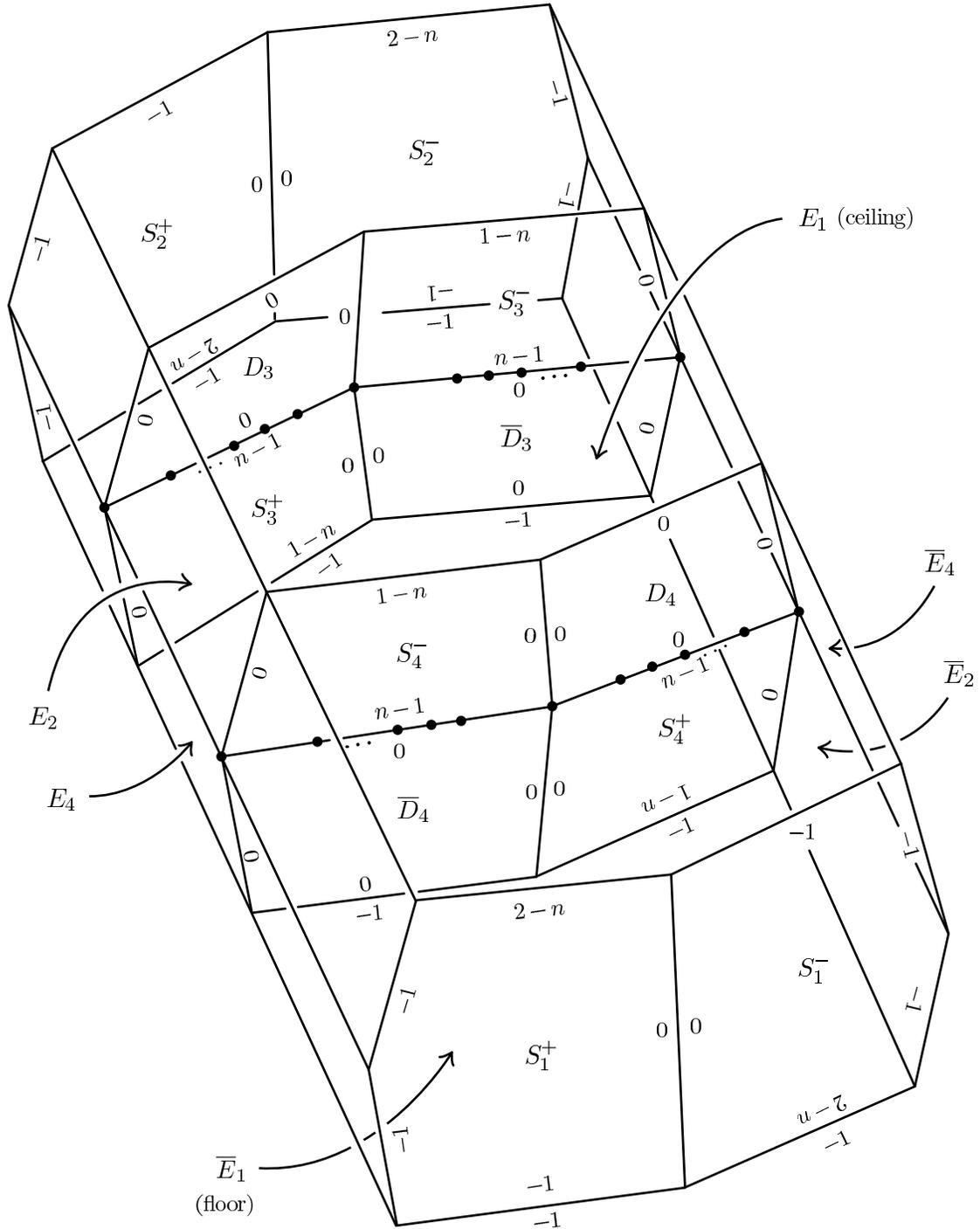}
\caption{structure of the conic bundle $X$ over $\tilde{\mathscr T}$ (Viewing this from above gives the projection to $\tilde{\mathscr T}$.)}
\label{fig_global}
\end{figure}

We note that the equation \eqref{c-bdle} makes sense globally over $\tilde{\mathscr T}$:
 one can readily verify 
 \begin{equation}\label{isom10}
 (\mathscr N\otimes\mathscr N')^{\vee}\simeq -(n-1)K_{\tilde{\mathscr T}}.
 \end{equation}
 On the other hand, there are $(n-1)$ anticanonical classes on the right hand side of \eqref{c-bdle}.
 Hence both sides belong to $H^0( -(n-1)K_{\tilde{\mathscr T}})$.
 In the following $X$ denotes the conic bundle defined by the equation \eqref{c-bdle}.
 (See Figure \ref{fig_global} which illustrates $X$ with some $\mathbf C^*$-invariant divisors.)
 
 \bigskip
\noindent
Proof of Theorem \ref{thm-1}.
As is already seen we have a surjective morphism  $\Phi_3:Z_3\to\tilde{\mathscr T}$ whose general fibers are smooth rational curves, equipped with two distinguished sections $E_1$ and $\ol{E}_1$.
We consider the direct image sheaf
\begin{align}
\mathscr E_n:=(\Phi_3)_*\mathscr O(E_1+\ol{E}_1).
\end{align}
By taking the direct image of the exact sequence
\begin{align}
0\,\lra\,\mathscr O\,\lra\,\mathscr O(E_1+\ol{E}_1)\,\lra\,
N_{E_1/Z_3}\oplus N_{\ol{E}_1/Z_3}\,\lra\,0,
\end{align}
we obtain the following exact sequence of sheaves on $\tilde{\mathscr T}$:
\begin{align}\label{di}
0\,\lra\,\mathscr O\,\lra\,\mathscr E_n\,\lra\,
\mathscr N\oplus \mathscr N'\,\lra\,R^1\Phi_3\,_*\mathscr O.
\end{align}
But we have $R^1\Phi_3\,_*\mathscr O=0$ since fibers of $\Phi_3$ are at most a string of rational curves.
Then the short exact sequence \eqref{di}  is readily seen to split, thanks to  \eqref{nb-1} and \eqref{nb-2}.
Hence we obtain
\begin{align}
\mathscr E_n\simeq\mathscr N\oplus\mathscr N'\oplus\mathscr O.
\end{align}
In particular, $\mathscr E_n$ is locally free.
Let $\mu:Z_3\to\mathbf P(\mathscr E_n^{\vee})$ be the meromorphic map associated to the pair of the morphism $\Phi_3$ and the line bundle $\mathscr O(E_1+\ol{E}_1)$.
On smooth fibers of $\Phi_3$, $\mu$ coincides with the rational map associated to the restriction of $\mathscr O(E_1+\ol{E}_1)$ to the fibers.
Hence  smooth fibers of $\Phi_3$ are mapped isomorphically to a conic in the fibers of $\mathbf P(\mathscr E_n^{\vee})$.
This means that $\mu$ is bimeromorphic.

Since  $\Phi_3$ and 
 $\mathscr O(E_1+\ol{E}_1)$ are $\mathbf C^*$-equivariant,
the target space $\mathbf P(\mathscr E_n^{\vee})$ of $\mu$ has a natural $\mathbf C^*$-action.
Recall that points of $E_1$  are $\mathbf C^*$-fixed.
So $\mathbf C^*$ acts fibers of the line bundle $N_{E_1/Z_3}$ by weight $1$ or $-1$, since otherwise the $\mathbf C^*$-action on $Z$ becomes non-effective or trivial.
If $\mathbf C^*$ acts on fibers of $N_{E_1/Z_3}$ by weight $1$ (resp.\,$-1$), it acts on fibers of  $N_{\ol{E}_1/Z_3}$ by weight $-1$ (resp.\,$1$) by reality.
(This can also be seen by explicitly calculating the action in a neighborhood of points of $C_1$ and $\ol{C}_1$ in the original twistor space $Z$.)
Hence we can suppose that the induced $\mathbf C^*$-action on $\mathscr E_n^{\vee}$ is of the form $(x,y,t)\mapsto (sx,s^{-1}y,t)$ for $s\in\mathbf C^*$, where $x$, $y$ and $t$ represent points of 
$\mathscr N^{\vee}$, $\mathscr N'^{\vee}$ and $\mathscr O$ respectively as in the proposition.
Since the image $\mu(Z_3)$ is $\mathbf C^*$-invariant in $\mathbf P(\mathscr E_n^{\vee})$,
this means that the defining equation of $\mu(Z_3)$  is of the form
\begin{align}
xy=R\,t^2
\end{align}
where $R$ is necessarily a section of $-(n-1)K_{\tilde{\mathscr T}}$ by \eqref{isom10}.
Further, since  the discriminant locus of $\Phi_3$ contains the anticanonical curve $\mathscr C_0$, $R$ can be divided by $P_0$ (= a defining equation of $\mathscr C_0$).

We have to show that the zero locus of $R/P_0\in H^0(-(n-2)K_{\tilde{\mathscr T}})$ decomposes into $(n-2)$ anticanonical curves.
To see this, we consider $\mathbf C^*$-fixed twistor lines in the twistor space $Z$.
As in Prop.\,\ref{prop-U(1)}, there are precisely $(n-2)$ twistor lines in $Z$ which have the property that all points are $\mathbf C^*$-fixed.
Let $L_j$ ($5\le j\le n+2$) be these fixed twistor lines.
(Note that $L_i=S_i^+\cap S_i^-$ $(1\le i\le 4)$ are other twistor lines.)
It is easy to see that these are disjoint from the cycle $C$.
In particular, the anticanonical map $\Phi:Z\to\mathscr T$ is a morphism on a neighborhood of $L_j$ ($5\le j\le n+2$).
Then the image $\Phi(L_j)$ cannot be a point since $\mathbf C^*$ acts non-trivially on every fiber of $\Phi$. 
Hence the image $\Phi(L_j)$ ($5\le j\le n+2$)
are curves in $\mathscr T$.
Therefore 
$\mathscr C_j:=\Phi_3(L_j)$ ($5\le j\le n+2$) are curves in $\tilde{\mathscr T}$
as well.
We show  that the surfaces $Y_j:=\Phi_3^{-1}(\mathscr C_j)$ ($5\le j\le n+2$) have singularities along $L_j$.
(Later it turns out each $Y_j$ consists of 2 irreducible components intersecting along $L_j$ transversally.)
To see this, we note that $Y_j$ is a $\mathbf C^*$-invariant divisor containing $L_j$.
On the other hand, from the fact that $L_j$ is a fixed twistor line, we can readily deduce that the natural $U(1)$-action is of the form $(u,v,w)\mapsto (su,s^{-1}v,w)$ $(s\in U(1)$) in a neighborhood of $L_j$, where $L_j$ is locally defined by $u=v=0$.
Hence any $\mathbf C^*$-invariant divisor containing $L_j$ must contain  (locally defined) divisors  $\{u=0\}$ or $\{v=0\}$.
In particular, our invariant divisor $Y_j$ contains at least one of these.
Moreover, since the map $\Phi_3$ is continuous in a neighborhood of $L_j$,
the images of these two divisors must be the same curve in $\tilde{\mathscr T}$.
This means that the surface $Y_j=\Phi_3^{-1}(\Phi_3(L_j))$ contains both of the two divisors.
Thus we see that $Y_j$ ($5\le j\le n+2$) have ordinary double points along general points of $L_j$, and any $\mathscr C_j$ are contained in the discriminant locus of $\Phi_3$.

Next we show that $\mathscr C_j$ $(5\le j\le n+2)$ are irreducible anticanonical curves on $\tilde{\mathscr T}$ which have a  unique node respectively.
Irreducibility is obvious since they are the images of the $\mathbf C^*$-fixed twistor lines.
To see that  $\mathscr C_j$ are anticanonical curves,  we first note that the anticanonical class $-K_{\mathscr T}$ of the original surface $\mathscr T$ is given by $\mathscr O_{\mathbf{CP}^4}(1)|_{\mathscr T}$, where $\mathscr T$ is embedded in $\mathbf{CP}^4$ as a quartic surface as before.
Therefore since $L_j\cdot(-K_Z)=4$, we obtain that $\Phi(L_j)\cdot(-K_{\mathscr T})$ is either $1,2$ or $4$ depending on the degree of the restriction  $L_j\to \mathscr \Phi(L_j)$ of $\Phi$.
Since the resolution $\tilde{\mathscr T}\to\mathscr T$ is crepant, the same is true for the intersection numbers $\mathscr C_j\cdot(-K_{\tilde{\mathscr T}})$.
On the other hand since $L_j$ $(5\le j\le n+2)$ intersect $S_i^+$ and $S_i^-$ ($1\le i\le 4)$, and since the intersection points are not on the cycle $C$, it follows that $\mathscr C_j$ actually intersects the image curves $\Phi_3(S_i^+)$ and  $\Phi_3(S_i^-)$ for any $1\le i\le 4$. 
As is already remarked, $\{\Phi_3(S_i^+),\Phi_3(S_i^-)\set 1\le i\le 4\}$ is precisely the set of irreducible components of reducible fibers for the projection $\tilde{\pi}:\tilde{\mathscr T}\to\Lambda\simeq\mathbf{CP}^1$.
Moreover since $\mathscr C_j\cdot(-K_{\tilde{\mathscr T}})\le 4$, the intersections $\Phi_3(S_i^+)\cap\mathscr C_j$ and $\Phi_3(S_i^-)\cap\mathscr C_j$ must be transversal (otherwise  $\mathscr C_j\cdot(-K_{\tilde{\mathscr T}})> 4$).
Thus we have seen that $\mathscr C_j$ ($5\le j\le n+2$) intersect $\Phi_3(S_i^+)$ and $\Phi_3(S_i^-)$ $(1\le i\le 4)$ transversally at a unique point respectively.
From this it readily follows that $\mathscr C_j$ are anticanonical curves of $\tilde{\mathscr T}$.
Then since it is an image of $L_j\simeq\mathbf{CP}^1$,  $\mathscr C_j$ has a unique singularity which is a node or a cusp.
But if it were a cusp, its inverse image by $\Phi_3$ would  be a real point.
Hence the singularity must be a node.

Thus we have shown that the discriminant locus of the morphism $\Phi_3$ contains the $(n-2)$ curves $\mathscr C_j\,(5\le j\le n+2)$ which are real anticanonical curves with a unique node respectively.
This means that $R/P_0$ can be divided by the product $P_5 P_6\cdots P_{n+2}$, where $P_j=0$ is a defining equation of $\mathscr C_j$. 
But since both are sections of $-(n-2)K_{\tilde{\mathscr T}}$ we can suppose $R/P_0=P_5 P_6\cdots P_{n+2}$.
Thus we have finished a proof of Theorem \ref{thm-1}.
\proofend

\section{Explicit construction of the twistor spaces}
In the last section we gave projective models of our twistor spaces as  conic bundles (Theorem \ref{thm-1}).
 In this section, reversing the procedures, we give an explicit construction of the twistor spaces starting from the projective models.
This is partially done already since we have given an explicit procedure for removing the indeterminacy of the anticanonical map $\Phi:Z\to\mathscr T\subset\mathbf{CP}^4$ and consequently obtained a morphism $\Phi_3:Z_3\to\tilde{\mathscr T}$ which is bimeromorphic to a conic bundle presented in Theorem \ref{thm-1}.
So it remains to analyze the bimeromorphic map $\mu:Z_3\to X\subset\mathbf P(\mathscr E_n^{\vee})$ obtained in the proof of Theorem \ref{thm-1}.
Recall that $\mu$ is the associated map to the pair $\Phi_3:Z_3\to\tilde{\mathscr T}$ and $\mathscr O(E_1+\ol{E}_1)$, so that there is a commutative diagram
\begin{equation}\label{cd4}
 \CD
Z_3@>{\mu}>> X \\
 @V{\Phi_3}VV @VV p V\\
\tilde{\mathscr T}@=\tilde{\mathscr T},\\
 \endCD
 \end{equation}
 where $p$ is the restriction of the projection $\mathbf P(\mathscr E_n^{\vee})\to\tilde{\mathscr T}$.
Note that $\Phi_3$ and $p$ are morphisms, but $\mu$ is a priori just a meromorphic map.

\begin{prop}\label{prop-alpha}
The bimeromorphic map $\mu$ satisfies the following.
(i) $\mu$ is a morphism.
(ii) $\mu$ contracts the 2 divisors $E_3$ and $\ol{E}_3$ in $Z_3$ to curves in $X$.
(iii) If $A$ is an irreducible divisor in $Z_3$ for which $\mu(A)$ is either a curve or a point, $A=E_3$ or $A=\ol{E}_3$ holds.
(iv) If $Z_3\to Z_4$ denotes the contraction of $E_3\cup E_3$ to curves as in (ii), then the induced morphism $Z_4\to X$ is a small resolution of all singularities of $X$.
\end{prop}

\noindent
Proof.
%First we show that the morphism $\Phi_4$ is equi-dimensional.
If $A$ is a divisor in $Z_3$ which is contracted to a point by the morphism $\Phi_3$, then its image into $Z$ (by the bimeromorphic morphism $Z_3\to Z$) becomes a divisor whose intersection number with twistor lines is zero.
Since such a divisor does not exist, $\Phi_3$ does not contract any divisor to a point.
Hence any fiber of $\Phi_3$ does not contain a divisor and therefore $\Phi_3$ is equi-dimensional.
Then since both $Z_3$ and $\tilde{\mathscr T}$ are non-singular,  $\Phi_3$ is a flat morphism.
This means that for any $y\in\tilde{\mathscr T}$, the restriction $\mu|_{\Phi_3^{-1}(y)}$ is precisely the rational map associated to the linear system
$\mathscr O(E_1+\ol{E}_1)|_{\Phi_3^{-1}(y)}$ (cf.\,\cite[pp.\,20--21]{U}).
Since a fiber $\Phi_3^{-1}(y)$ is at most a chain of rational curves, this implies that 
$\mu|_{\Phi_3^{-1}(y)}$ is a morphism contracting components which do not intersect the sections $E_1$ nor $\ol{E}_1$.
This in particular means $\mu$ is a morphism and we obtain (i).
Also since the divisors $E_3$ and $\ol{E}_3$ are disjoint form $E_1\cup\ol{E}_1$, it follows that $E_3$ and $\ol{E}_3$ are contracted to curves by $\mu$.
Hence we obtain (ii).
It is readily seen by our explicitness of the bimeromorphic morphism $\mu_3:Z_3\to Z$,  that $E_3$ and $\ol{E}_3$ in $Z_3$ are  biholomorphic to $\mathbf{CP}^1\times\mathbf{CP}^1$ and that their normal bundles in $Z_3$ have degree $(-1)$ along  directions of the contractions (cf.\,(c) of Figures \ref{fig1_and_2},\ref{fig3},\ref{fig4}).
Therefore the morphism $\mu$ blows down $E_3$ and $\ol{E}_3$ to curves in $X$
in a way that  the resulting 3-fold $Z_4$ is still non-singular
((c) $\to$ (d) of Figures \ref{fig1_and_2},\ref{fig3},\ref{fig4}).
Since $E_3$ and $\ol{E}_3$ were originally exceptional divisors over $C_3$ and $\ol{C}_3$, and $\Phi(C_3)$ and $\Phi(\ol{C}_3)$ are the nodes of the surface $\mathscr T$, 
the curves $\Phi_3(E_3)$ and $\Phi_3(\ol{E}_3)$ must be the exceptional curves $\Gamma$ and $\ol{\Gamma}$ of the resolution $\tilde{\mathscr T}\to\mathscr T$.
(See the diagram \eqref{string}.)

Next to show (iii)
suppose that $A$ is an irreducible divisor in $Z_3$ which is contracted  by $\mu$
and that $A$ is different from $E_3$ and $\ol{E}_3$.
Then since $\mu$ is $\mathbf C^*$-equivariant and its image $X$ is 3-dimensional, 
$A$ must be $\mathbf C^*$-invariant.
Further $A$ cannot be an irreducible component of  the bimeromorphic morphism
$\mu_3:Z_3\to Z$, since all of them intersect at least one of $E_1$ and $\ol{E}_1$ along a curve, except $E_3$ and $\ol{E}_3$, so that they cannot be contracted by $\mu$.
Hence the image $\mu_3(A)$  must be a $\mathbf C^*$-invariant divisor in $Z$.
Further it cannot be an irreducible component of a member of the pencil $|(-1/2)K_Z|$, since such a component always contains at least one of $C_1$ and $\ol{C}_1$ and hence $A$ intersects at least one of $E_1$ and $\ol{E}_1$ along a curve.
This means that  the image $\Phi_3(A)$ is a curve intersecting any fibers of $\tilde{\pi}:\tilde{\mathscr T}\to\Lambda$.
Take a general real member $S$ of $|(-1/2)K_Z|$ and let $S'$ be its strict transform in $Z_3$.
$S'$ is biholomorphic to $S$.
Then $\mu|_{S'}$ is precisely the contraction of the intersection curves
$S'\cap E_3\,(\simeq C_3\subset S)$ and $S'\cap\ol{E}_3\,(\simeq \ol{C}_3\subset S)$, since on $S$ there is no $\mathbf C^*$-orbit which does not intersect $C_1$ nor $\ol{C}_1$, except $C_3$ and $\ol{C}_3$.
Because $S$ can be supposed to be a general member, this implies that the curve $\Phi_3(A)$ must coincide with $\Gamma$ or $\ol{\Gamma}$.
Hence by \eqref{dl3}, $A$ must be  one of $E_i$ and $\ol{E}_i$, $2\le i\le 4$.
But these cannot happen.
Thus we obtain (iii).

The final assertion (iv) is obvious since $\mu$ does not contract any divisor and the 3-fold $Z_4$ obtained by contracting $E_3$ and $\ol{E}_3$ is non-singular.
\proofend

\bigskip
By Prop.\,\ref{prop-alpha} and its proof, the bimeromorphic map $\mu$ is the composition of the blowing-down of the divisors $E_3$ and $\ol{E}_3$ 
along the directions of the projections to $\Gamma$ and $\ol{\Gamma}$ respectively,
with contractions of some $\mathbf C^*$-invariant rational curves (which are necessarily disjoint from $E_1$ and $\ol{E}_1$)  into isolated singularities of $X$.
In the equation \eqref{c-bdle} of $X$ the image curves $\mu(E_3)$ and $\mu(\ol{E}_3)$ are contained in the reducible curve $\{x=y=P_0=0\}$ (which are mapped biholomorphically to the curve $\mathscr C_0$ by $p:X\to\tilde{\mathscr T}$), and the singularities  of $X$ are over the singularities of the discriminant locus of $p$.
Of course, the latter singularities are either the intersection points of the anticanonical curves $\mathscr C_j=\{P_j=0\}$, where $j=0$ or $5\le j\le n+2$, or the singularities of $\mathscr C_j$ themselves.
$\mathscr C_0$ is a cycle of 6 rational curves.
Hence it has 6 ordinary nodes.
On the other hand $\mathscr C_j$ $(5\le j\le n+2)$ have a unique node respectively.
Thus in general the number of the  singularities of $X$ is
\begin{align}
6+(n-2)+4
\begin{pmatrix}
n-1\\ 2
\end{pmatrix}.
\end{align}
(Note that $\mathscr C_j\cdot\mathscr C_k=(-K_{\tilde{\mathscr T}})^2=4$ which results in the last term.)

Reversing all the operations we have obtained, it is  now possible to give explicit way for obtaining our twistor space $Z$ from the projective model $X$.
It can be summarized in the following diagram
\begin{equation}\label{string2}
 \CD
Z_3@>>> Z_2@>>>  Z_1@>>> Z\\
   @VVV             @.         @VV\Phi_1V @VV\Phi V\\
 Z_4 @>>> X@>>{p}>              \tilde{\mathscr T}@>>>\mathscr T.\\
 \endCD
 \end{equation}
% Here we briefly describe how to obtain our twistor space $Z$.
  Namely all the operations can be briefly described as follows:
 \begin{itemize}
 \item
 Suppose $n\ge 4$ and let $Z$ be a twistor space on $n\mathbf{CP}^2$ as in Theorem \ref{thm-1}.
 We want an explicit construction of $Z$.
 (The condition $n\ge 4$ is superfluous, and the following construction perfectly works for $n\ge 2$. See   final part of \S 5.2.) 
 \item
We fix an integer $n\ge 4$ and let $\alpha$ be a real number satisfying $-1<\alpha<0$.
Then the quartic surface $\mathscr T$ defined by the equation \eqref{T} is determined.
($\mathscr T$ is independent of $n$.)
$\mathscr T$ will serve as a minitwistor space.
Let $\tilde{\mathscr T}\to\mathscr T$ be the minimal resolution of the 2 nodes of $\mathscr T$ (cf.\,Prop.\,\ref{prop-T}).
\item
Next we realize $\tilde{\mathscr T}$ as 4 points blown-up of $\mathbf{CP}^1\times\mathbf{CP}^1$ by the blowing-down $\nu:\tilde{\mathscr T}\to\mathbf{CP}^1\times\mathbf{CP}^1$ explicitly given for obtaining \eqref{nb-1}.
Then we can consider the line bundles $\mathscr N=\nu^*\mathscr O(-1,2-n)$ and $\mathscr N'=\ol{\sigma^*\mathscr N}$ over $\tilde{\mathscr T}$, where $\sigma$ denotes the natural real structure on $\tilde{\mathscr T}$.
(In the homogeneous coordinate  of Prop.\,\ref{prop-quotient},
$\sigma$ is explicitly given by $(y_0,y_1,y_2,y_3,y_4)\mapsto (\ol{y}_0,\ol{y}_1,\ol{y}_2,\ol{y}_4,\ol{y}_3)$.)
We put $\mathscr E_n=\mathscr N\oplus\mathscr N'\oplus \mathscr O$, a rank-3 vector bundle on $\tilde{\mathscr T}$.
\item
Let $X$ be a conic bundle in $\mathbf P(\mathscr E_n^{\vee})$ defined by the equation \eqref{c-bdle}, where $P_0$ and $P_j$ ($5\le j\le n+2)$ are real sections of the anticanonical bundle $-K_{\tilde{\mathscr T}}$ such that the zero locus $\{P_0=0\}$ is the anticanonical  curve $\mathscr C_0$ introduced in the explanation before Theorem \ref{thm-1}
(see also Figure \ref{figT}), and that $\mathscr C_j=\{P_j=0\}$ ($5\le j\le n+2)$ are real irreducible curves with a unique node respectively. 
(See Figure \ref{fig_global}.)
Further we suppose that all the intersections of $\mathscr C_j$'s ($j=0$ or $5\le j\le n+2$) are transversal.
Let $p:X\to\tilde{\mathscr T}$ be the natural projection.
All the singularities of $X$ are ordinary double points  over the singularities of the curve $\mathscr C_0\cup\mathscr C_5\cup\cdots\cup\mathscr C_{n+2}$.
They are also lying on the section $\{x=y=0\}$ of the bundle  $\mathbf P(\mathscr E_n^{\vee})\to\tilde{\mathscr T}$.
(In Figure \ref{fig_global}, all ODP's lying over the curve $\mathscr C_0$ are denoted by dotted points.)
\item
Let $Z_4\to X$ be  small resolutions of all the ordinary double points of $X$.
For the ODP's on $S_3^{\pm}$ and $S_4^{\pm}$, we choose small resolutions which blow-up $S_3^{\pm}$ and $S_4^{\pm}$.)
Let $Z_3\to Z_4$ be a blowing-up along two smooth rational curves 
$E_2\cap E_4$ and $\ol{E}_2\cap \ol{E}_4$ (see Figure \ref{fig_global}).
%determined by the following condition: their images in $X$ are smooth rational curves which are mapped isomorphically to the curves $\Gamma$ or $\ol{\Gamma}$, where  $\Gamma$ and $\ol{\Gamma}$ are the exceptional curves of the resolution $\tilde{\mathscr T}\to\mathscr T$;
%further, the image curves in $X$ are on the section $\{x=y=0\}$.
\item
Let $Z_3\to Z_2$ be the blowing-down of the 4 divisors $D_3,\ol{D}_3,D_4,\ol{D}_4$ to curves, along the directions which accord  the projection $\tilde{\pi}:\tilde{\mathscr T}\to\Lambda\simeq\mathbf{CP}^1$
(cf.\,(c)$\to$(b) of Figures \ref{fig3},\ref{fig4}).
$Z_2$ is still non-singular.
\item 
In $Z_2$, the divisors $S_i^{\pm}$ ($1\le i\le 4$)  respectively possesses two $(-1,-1)$-curves as in (b) of Figures \ref{fig1_and_2},\ref{fig3},\ref{fig4}. Let $Z_2\to Z_1$ be the contraction of the curves into ordinary double points.
Then we are in the situation that  the divisor $\sum_{i=1}^4(E_i+\ol{E}_i)$ can be simultaneously blown-down to a cycle $C=\sum_{i=1}^4(C_i+\ol{C}_i)$ of rational curves.
Let $Z_1\to Z$ be the blowing-down.
Then $Z$ is the space we are seeking.
 \end{itemize}

\section{Existence, moduli, and similarities with LeBrun twistor spaces}
\noindent
{\bf (5.1)}
In this subsection we first show that our twistor spaces can be obtained as a $\mathbf C^*$-equivariant deformation of the twistor space of some Joyce metric \cite{J95} on $n\mathbf{CP}^2$.
Next we show that the property of our twistor spaces that they possess a rational surface $S$ constructed by Section 2 as a member of the system $|(-1/2)K_Z|$, is preserved under $\mathbf C^*$-equivariant small deformations.
Since all our results rely on the existence of this divisor, it means that the structure of our twistor spaces is stable under $\mathbf C^*$-equivariant small deformations.

First we explain which Joyce metric we shall consider.
For this it suffices to specify the structure of a smooth toric surface which is contained in the twistor space as a torus invariant member of the system $|(-1/2)K_Z|$.
For constructing the toric surface explicitly, as in the construction of our surface $S$ in Section 2, we start from $\mathbf{CP}^1\times\mathbf{CP}^1$ and choose a non-real member $C_1$ of $|\mathscr O(1,0)|$.
Next we choose a point $p_1\in C_1$ and we blow-up $\mathbf{CP}^1\times\mathbf{CP}^1$ at $p_1$ and $\ol{p}_1$ (where $\ol{p}_1$ is the image of $p_1$ under the real structure \eqref{rs5}).
Next we blow-up the resulting surface at the two intersection points of $C_1\cup \ol{C}_1$ and the exceptional curves.
Repeating this blowing-up procedure $(n-1)$ times, we obtain a toric surface with $c_1^2=8-2(n-1)$,
where $C_1$ and $\ol{C}_1$ satisfy $C_1^2=\ol{C}_1^2=1-n$ on the surface. 
The exceptional curves of the birational morphism onto the original surface $\mathbf{CP}^1\times\mathbf{CP}^1$ contains precisely four $(-1)$-curves, and all of them intersect $C_1\cup \ol{C}_1$.
As the final step for obtaining the toric surface, we choose a conjugate pair of $(-1)$-curves among these, and blow-up the torus-invariant points on the curves which do not on $C_1\cup\ol{C}_1$.pro
Let $S_J$ be a toric surface obtained in this way.
$S_J$ satisfies $c_1^2=8-2n$.

Then it is elementary to see that our rational surface $S$ with $\mathbf C^*$-action (explicitly constructed in Section 2) is obtained as a $\mathbf C^*$-equivariant deformation 
of the toric surface $S_J$, where $\mathbf C^*$-subgroup of $(\mathbf C^*)^2$ of the toric surface is the one specified by the condition that it fixes points of $C_1$ and $\ol{C}_1$.

On the other hand, let $Z_J$ be the twistor space of a Joyce metric on $n\mathbf{CP}^2$ 
which has the toric surface $S_J$ as a ($(\mathbf C^*)^2$-invariant) member of $|(-1/2)K|$.
Then the following result means that our twistor spaces studied in Section 2--4 can be obtained as a $\mathbf C^*$-equivariant small deformation  of $Z_J$,
where $\mathbf C^*\subset(\mathbf C^*)^2$ is the subgroup chosen in the last paragraph.

\begin{prop}
Let $Z_J$ be the twistor space of a Joyce metric containing the toric surface $S_J$ as above, and $\mathbf C^*\subset(\mathbf C^*)^2$ the subgroup specified as above.
 Then if $Z$ is a twistor space obtained as a small $\mathbf C^*$-equivariant deformation of $Z_J$, then $Z$ contains the surface $S$ as its member of the system $|(-1/2)K_Z|$.
\end{prop}

\noindent Proof.
Since the twistor space of a Joyce metric is Moishezon, we have $H^2(\Theta_{Z_J}\otimes\mathscr O(-S_J))=0$ \cite[Lemma 1.9]{C91}.
By a result of Horikawa \cite{Hor76} this implies that the surface $S_J$ is costable under small  deformations of $Z_J$.
Namely for any small deformations of $S_J$, there exists a deformation of  the pair $(Z_J,S_J)$ such that deformation of $S_J$ coincides with the given one.
This is also true for equivariant deformations.
By applying this to the above deformation of $S_J$ into $S$, we obtain the required twistor space $Z$ having $S$ as a member of $|(-1/2)K_Z|$.
\proofend

\bigskip
Next we show that the structure of our twistor space $Z$ is stable under $\mathbf C^*$-equivariant small deformations.
Namely, we show that the structure of the member $S\in |(-1/2)K_Z|$ is stable under any $\mathbf C^*$-equivariant small deformations of $S$, and $S$ always survives under $\mathbf C^*$-equivariant small deformations of the twistor space $Z$.

It is generally true that if $S$ is a complex surface equipped with $\mathbf C^*$-action and if $S$ is obtained from another rigid complex surface $S_1$ with $\mathbf C^*$-action by a succession of blowing-up at $\mathbf C^*$-fixed points, then any $\mathbf C^*$-equivariant small deformations of $S$ are obtained as a deformation obtained by moving the  blown-up points on $S_1$.
(This  is due to `the stability of a $(-1)$-curve' under small deformations of a surface \cite{Kod}, plus the supposed rigidity of the starting surface $S_1$.)
For our complex surface $S$, among $2n$ points of $\mathbf{CP}^1\times\mathbf{CP}^1$ to be blown-up, the two points of the final blowup $S\to S'$ cannot be moved since they are isolated $\mathbf C^*$-fixed points.
The remaining $2(n-1)$ points are on $C_1$ and $\ol{C}_1$ which are $\mathbf C^*$-fixed locus.
Hence any $\mathbf C^*$-equivariant small deformations of our rational surface $S$ are  obtained by moving these $2(n-1)$-points on $C_1\cup\ol{C}_1$.
Thus the structure of $S$ does not change while we are considering $\mathbf C^*$-equivariant deformations.

Next we see that $S$ survives under $\mathbf C^*$-equivariant small deformations of $Z$.
By an equivariant version of a criterion of Kodaira \cite{Kod} about stability of submanifolds under small deformations of  ambient space, 
it suffices to verify that our rational surface $S$ satisfies $H^1(-K_S)^{\mathbf C^*}=0$, since we have $N_{S/Z}\simeq -K_S$.
For this, we have an obvious exact sequence 
\begin{align}
0\,\lra\,-K_S-C_1-\ol{C}_1\,\lra\,-K_S\,\lra\, -K_S|_{C_1+\ol{C}_1}\,\lra\,0,
\end{align}
and we have $-K_S|_{C_1+\ol{C}_1}\simeq\mathscr O_{C_1}(3-n)\oplus\mathscr O_{\ol{C}_1}(3-n)$ by adjunction formula.
It is also  routine computations to see that $H^i(-K_S-C_1-\ol{C}_1)=0$ for $i=1,2$.
Hence we obtain
an equivariant isomorphism 
\begin{align}
H^1(-K_S)\simeq H^1(-K_S|_{C_1+\ol{C}_1}).
%\,\,\,(\simeq \mathbf C^{n-4}\oplus\mathbf C^{n-4})
\end{align}
Further, by our explicit construction of $S$ it is possible to compute the natural $\mathbf C^*$-action on $H^1(-K_S|_{C_1})\simeq\mathbf C^{n-4}$ and  $H^1(-K_S|_{\ol{C}_1})\simeq \mathbf C^{n-4}$ concretely.
The result is that all  the weights of the $\mathbf C^*$-action on $H^1(-K_S|_{C_1})$ are either $1$ or $-1$,
and  those for the $\mathbf C^*$-action on $H^1(-K_S|_{\ol{C}_1})$ are $-1$ or $1$ respectively.
Hence there is no $\mathbf C^*$-invariant point on $H^1(-K_S)$ other than $0$.
Thus we obtain that $S$ survives under $\mathbf C^*$-equivariant small deformations of $Z$.
(Of course, the complex structure of $S$ itself deforms under $\mathbf C^*$-equivariant deformations, if we perturb $(n-1)$ chosen points along $C_1\in|\mathscr O(1,0)|$ to be blown-up.)

\bigskip
\noindent
{\bf (5.2)}
In this subsection we compute the dimension of the moduli space of our twistor spaces studied in Sections 2--4, by counting the number of parameters involved in our construction of the twistor spaces.
Also we see that (if $n\ge 4$) our twistor spaces cannot be obtained as a small $\mathbf C^*$-equivariant deformation of LeBrun twistor spaces of any kind.
Finally we give a remark for the case $n=2,3$.

Recall that the projective models of our twistor spaces have a structure of conic bundles over the rational surface $\tilde{\mathscr T}$, and that  $\tilde{\mathscr T}$ is the minimal resolution of  the quartic surface $\mathscr T$ 
 defined by the equations
\begin{equation}\label{T2}
y_1y_2=y_0^2,\,\,\,y_3y_4=y_0\{y_1-\alpha y_2+(\alpha-1 )y_0\},
\end{equation}
where $\alpha$ satisfies $-1<\alpha<0$ (Prop.\,\ref{prop-quotient}).
The complex structure of $\mathscr T$ deforms if we move $\alpha$, 
since $\alpha$ corresponds to one of the 4 discriminant points of the natural projection  $\tilde{\pi}:\tilde{\mathscr T}\to\Lambda$, and the remaining 3 discriminant points are fixed (cf. Prop.\,\ref{prop-T} and its proof).
Thus we have one parameter for specifying the base surface $\mathscr T$ or $\tilde{\mathscr T}$.

Next, once we fix $\alpha$, the projective models of the conic bundles are defined by the equation
\begin{align}\label{c-bdle2}
xy=t^2P_0P_5P_6\cdots P_{n+2},
\end{align}
where $P_0$ and $P_j$ ($5\le j\le n+2$) are anticanonical curves on the surface $\tilde{\mathscr T}$ (Theorem \ref{thm-1}).
In particular, the conic bundles are uniquely determined by the anticanonical curves
\begin{align}
\mathscr C_j=\{P_j=0\},\,\,\,j=0\,\,\,{\text{or}}\,\,5\le j\le n+2.
\end{align}
Among these $(n-1)$ curves, $\mathscr C_0$ is a cycle of 6 rational curves, and it is  uniquely determined from the complex surface $\tilde{\mathscr T}$.
(Namely its irreducible components consist of the exceptional curves of the resolution $\tilde{\mathscr T}\to\mathscr T$ and two of the reducible fibers of the projection $\tilde{\pi}:\tilde{\mathscr T}\to\mathscr T$. See Figure \ref{fig_global}).
So there is no freedom in determining $\mathscr C_0$.
On the other hand, for $\mathscr C_j$ ($5\le j\le n+2$), 
we note that $\dim |-K_{\tilde{\mathscr T}}|=4$ since $-K_{\mathscr T}\simeq\mathscr O_{\mathbf{CP}^4}(1)|_{\mathscr T}$ and the resolution $\tilde{\mathscr T}\to\mathscr T$ is crepant.
However, $\mathscr C_j$ ($5\le j\le n+2$) are not general anticanonical curves but 
have a unique node respectively.
This drops 1-dimension and for each $j$ ($5\le j\le n+2$)  there are 3-dimensional freedom  of choices.
Thus the number of parameters for fixing discriminant locus is 
$3(n-2)$.
Further, the identity component for the group of holomorphic  automorphism of $\mathscr T$ or $\tilde{\mathscr T}$ is $\mathbf C^*$, where it explicitly acts by
\begin{align}
(y_0,y_1,y_2,y_3,y_4)\longmapsto (y_0,y_1,y_2,sy_3,s^{-1}y_4),\,\,\,s\in\mathbf C^*
\end{align} 
in the coordinate \eqref{T2} on $\mathbf{CP}^4$.
This $\mathbf C^*$-action naturally induces that on the space of the choices of $\mathscr C_j$ ($5\le j\le n+2$).
Summing these up, the dimension of the moduli space of our twistor space is 
\begin{align}\label{dim}
\{1+3(n-2)\}-1=3n-6.
\end{align}
Note that this is identical to the dimension of the moduli space of LeBrun metrics (on $n\mathbf{CP}^2$) whose identity component of the automorphism group is precisely $U(1)$.

Next we remark that our twistor spaces on $n\mathbf{CP}^2$ cannot be obtained as a $\mathbf C^*$-equivariant small deformation of LeBrun twistor spaces. 
For general LeBrun twistor spaces on which only $\mathbf C^*$ acts effectively, 
this is actually true since equivariant deformations of such LeBrun twistor spaces are still LeBrun twistor spaces \cite{LB93}.
For LeBrun twistor spaces admitting an effective $(\mathbf C^*)^2$-action, 
it is determined in \cite{Hon07-1} which  $\mathbf C^*$-subgroup  (of $(\mathbf C^*)^2$) admits equivariant deformation such that the resulting space is not a LeBrun twistor space.
It is also shown that  the moduli spaces of twistor spaces on $n\mathbf{CP}^2$ obtained this way is either $n$ or $n+2$.
If $n\ge 5$, these cannot be equal to the dimension \eqref{dim}.
Hence
 our twistor spaces cannot be obtained as a $\mathbf C^*$-equivariant deformation of LeBrun metrics of any kinds, for the case $n\ge 5$.
If $n=4$, we have $n+2=3n-6$ and the dimensions of the moduli spaces coincide.
But the results of  \cite{Hon01} and \cite{Hon07-1} show that small deformations of LeBrun metrics on $4\mathbf{CP}^2$ with $(\mathbf C^*)^2$-action which have 6-dimensional moduli
always drop algebraic dimension of the twistor spaces.
Since our twistor spaces are of course Moishezon, this means that the above conclusion (for $n\ge 5$) is true also for the case $n=4$.

So far in this paper we have always supposed  that the  twistor spaces are on $n\mathbf{CP}^2$, $n\ge 4$.
But the results of Sections 2--4 can be readily justified for the case $n=2$ and $n=3$ 
if we use the $\mathbf C^*$-fixed part $|-K_Z|^{\mathbf C^*}$ instead of $|-K_Z|$.
Consequently, the explicit  construction in Section 4 works also for the case $n=2,3$.
(If $n=2$, we read the equation \eqref{c-bdle} as  `$xy=t^2P_0$'.
Namely $\mathscr C_0$ becomes the unique discriminant anticanonical curve of the conic bundle and no nodal components appear.)
On the other hand, the construction does not work if $n=0$ or $n=1$, since in these cases the construction of our starting surface $S$ has no meaning.

If $n=2$, the twistor spaces we obtain are of course nothing but Poon's twistor spaces studied in \cite{P86}.
Thus for the case $n=2$  our construction gives a new realization of Poon's twistor spaces.
In \cite{P86} Poon used the system $|(-1/2)K_Z|$ and showed that it induces a bimeromorphic morphism from $Z$ to a quartic 3-fold in $\mathbf{CP}^5$.
On the other hand our study is based on the system $|-K_Z|^{\mathbf C^*}$ which yielded the meromorphic map $\Phi:Z\to \mathscr T\subset\mathbf{CP}^4$.

If $n=3$, since our twistor spaces admit $\mathbf C^*$-action, the twistor spaces are either LeBrun twistor spaces or the twistor spaces of double solid type studied in \cite{Hon07-2}.
But since  $|(-1/2)K_Z|^{\mathbf C^*}$ is only a pencil, they cannot be LeBrun twistor spaces. 
(For LeBrun twistor spaces the system is always 3-dimensional.)
Thus for the case $n=3$  the present construction yields another realization of the twistor spaces in \cite{Hon07-2}.

\bigskip
\noindent
{\bf (5.3)}
As one may notice, our twistor spaces resemble  LeBrun twistor spaces \cite{LB91} in many respects.
In this subsection we discuss  these similarities   in detail, as well as their differences.
Throughout this section  $Z_{\,\rm{LB}}$ denotes a LeBrun twistor space on $n\mathbf{CP}^2$.
Then the half-anticanonical system induces a meromorphic map
$$\Phi_{\,\rm{LB}}:Z_{\,\rm{LB}}\to\mathbf{CP}^3$$
whose image is a non-degenerate quadratic surface $\mathscr T_{\,\rm{LB}}\simeq \mathbf{CP}^1\times\mathbf{CP}^1$.
$Z_{\,\rm{LB}}$ also admits a $\mathbf C^*$-action and the map
$\Phi_{\,\rm{LB}}$ is $\mathbf C^*$-equivariant, where $\mathbf C^*$ acts trivially on the target space.
Further, general fibers of $\Phi_{\,\rm{LB}}$ are the closures of orbits which are irreducible smooth rational curves.
Thus our map $\Phi:Z\to\mathscr T\subset\mathbf{CP}^4$ can be thought as an analogue of $\Phi_{\,\rm{LB}}:Z_{\,\rm{LB}}\to\mathscr T_{\,\rm{LB}}\subset\mathbf{CP}^3$,
and the surface $\mathscr T$ is an analogue of $\mathscr T_{\,\rm{LB}}$.
%Hence like our meromorphic map $\Phi:Z\to\mathscr T$,  the map $\Phi_{\,\rm{LB}}:Z_{\,\rm{LB}}\to\mathscr T_{\,\rm{LB}}$ can be regarded as a  quotient map of the $\mathbf C^*$-action.
One of the differences is that while  $\Phi_{\,\rm{LB}}$ is the meromorphic map associated to the system $|(-1/2)K_Z|$, our map $\Phi$ is associated to the twice $|-K_Z|$.
Further,  $\mathscr T_{\,\rm{LB}}$ is smooth, while our surface $\mathscr T$ has 2 ordinary double points as in Prop.\,\ref{prop-T}.
More significantly, the defining equations \eqref{T} of  $\mathscr T$ contain a parameter $\alpha$ and as explained in \S 5.2 the complex structure of the surface actually deforms if we move the parameter.
Thus our image surface constitute a 1-dimensional moduli space.
(In \cite{Hon-p-1} it was shown that the moduli space can be identified with the moduli space of elliptic curves defined over real numbers.)
In contrast, for LeBrun twistor spaces the image surface $\mathscr T_{\,\rm{LB}}$ is rigid, of course.

Secondly, for LeBrun twistor spaces, the indeterminacy locus of $\Phi_{\,\rm{LB}}$ (that is,  the base locus of  $|(-1/2)K_Z|$) is a conjugate pair of smooth rational curves (which are contained in the fixed locus of the $\mathbf C^*$-action), and 
if we blow up these curves with the resulting space  $Z'_{\,\rm{LB}}$,
the map $\Phi_{\,\rm{LB}}$ already becomes a morphism  $Z'_{\,\rm{LB}}\to\mathscr T_{\,\rm{LB}} $ having the exceptional divisors as sections.
In contrast, the indeterminacy locus of our map $\Phi$ (i.\,e. the base locus of the anticanonical system) is somewhat complicated as in Cor.\,\ref{cor-ac}, and we had to make a succession of blow-ups to eliminate the indeterminacy, as in the sequence \eqref{string}.
This is the reason why we need some complicated construction, compared to the construction of LeBrun twistor spaces.
%Thus although the situation is pretty complicated, what we have done is 

Thirdly we explain similarities on defining equations of  projective models of the twistor spaces.
For LeBrun twistor spaces, there is a bimeromorphic morphism $\mu_{\,\text{LB}}$ from the blown-up space $Z'_{\,\rm{LB}}$   to the conic bundle
\begin{align}\label{LB}
X_{\,{\rm LB}}:\,xy=P_1P_2\cdots P_n \,t^2,
\end{align}
where $x\in\mathscr O(1,n-1)$, $y\in\mathscr O(n-1,1)$, $t\in\mathscr O$, over the  surface $\mathscr T_{\,\rm{LB}}\simeq\mathbf{CP}^1\times\mathbf{CP}^1$,
and $P_1,\cdots, P_n$ are  real sections of $\mathscr O(1,1)$.
This bimeromorphic morphism $\mu_{\,\text{LB}}$ is obtained \cite{Ku92} as the canonical map associated to the pair of the morphism $Z'_{\,\rm{LB}}\to\mathscr T_{\,\rm{LB}}$ and the line bundle $\mathscr O(E_1+\ol{E}_1)$, where $E_1$ and $E_1$ are sections of the morphism $Z'_{\,\rm{LB}}\to\mathscr T_{\,\rm{LB}}$ which are the exceptional divisors of the blowing-up $Z'_{\,\rm{LB}}\to Z_{\,\rm{LB}}$.
Thus our bimeromorphic morphism $\mu:Z_3\to X$ studied in the proof of Prop.\,\ref{prop-alpha} is an analogue of  $\mu_{\,\text{LB}}:Z'_{\,\rm{LB}}\to X_{\,\rm{LB}}$.
But note that while $\mu_{\,\text{LB}}$ is exactly small resolutions of singularities of $X_{\,\text{LB}}$,
our morphism $\mu$  contracts not only curves but also the 2 divisors $E_3$ and $\ol{E}_3$.

The equation \eqref{LB} also shows that  the discriminant locus of the conic bundle 
$X_{\,\rm{LB}}\to\mathscr T_{\,\rm{LB}}$ splits into $n$ irreducible components $\{P_j=0\}$, $1\le j\le n$.
This is very similar to our case (Theorem \ref{thm-1}).
Moreover, it
 is also true that the curves $\{P_j=0\}\subset\mathscr T_{\,\rm{LB}}$ $(1\le j\le n)$ are exactly the images of $\mathbf C^*$-fixed twistor lines.
As showed in the proof of Theorem \ref{thm-1} this is true also for our curves $\mathscr C_j=\{P_j=0\}$ for $5\le j\le n+2$.
(For LeBrun twistor spaces there are $n$ such twistor lines;
 in our case there are only $(n-2)$ such twistor lines as stated in Prop.\,\ref{prop-u1}.)
On the other hand big difference is that among anticanonical curves in the discriminant locus, one component $\{P_0=0\}$ plays a special role in our case;
namely although it is an anticanonical curve like other components, it consists of 6 irreducible components.
In LeBrun's case, there is no such special one among $P_j=0$.

Finally we mention a remarkable difference about the images of twistor lines.
For LeBrun twistor spaces, the image $\Phi_{\,{\rm LB}}(L)$ of a general twistor line $L$ is a non-singular curve (whose bidegree is $(1,1)$).
On the other hand we showed in the proof of Theorem \ref{thm-1} that $\Phi(L_j)$, $5\le j\le n+2$, are {\em nodal}\, rational curves.
The proof works not only for $L_j$ but also for generic twistor lines $L$.
Consequently $\Phi(L)$ is a nodal rational (and anticanonical) curve in $\mathscr T$.
(For another proof of this fact, see \cite{Hon-p-1}.)
In short, general minitwistor lines in our minitwistor space $\mathscr T$ are nodal rational curves in the anticanonical class.

\small
\vspace{13mm}
\hspace{7.5cm}
$\begin{array}{l}
\mbox{Department of Mathematics}\\
\mbox{Graduate School of Science and Engineering}\\
\mbox{Tokyo Institute of Technology}\\
\mbox{2-12-1, O-okayama, Meguro, 152-8551, JAPAN}\\
\mbox{{\tt {honda@math.titech.ac.jp}}}
\end{array}$

\end{document}